\documentclass[10pt]{amsart}

\usepackage{amssymb}
\usepackage{amsmath}
\usepackage{amsthm}
\usepackage{stmaryrd}
\usepackage{verbatim}
\usepackage{listings}
\usepackage[hidelinks]{hyperref}
\usepackage{enumerate}
\usepackage{array}
\usepackage{multirow}
\usepackage{nicefrac}

\usepackage{longtable}
\usepackage[backrefs,lite]{amsrefs}
\DefineSimpleKey{bib}{primaryclass}{}
\DefineSimpleKey{bib}{archiveprefix}{}

\BibSpec{arXiv}{%
+{}{\PrintAuthors}{author}
+{,}{ \textit}{title}
+{}{ \parenthesize}{date} +{,}{ arXiv }{eprint} +{,}{ primary class
}{primaryclass} }

\usepackage{tikz, pgfplots, pgfplotstable}
\usetikzlibrary{cd}
\usetikzlibrary{calc}
\usetikzlibrary{decorations.markings}
\usetikzlibrary{positioning,arrows}
\usepgfplotslibrary{external}
\usepackage{etoolbox}
\AtBeginEnvironment{tikzcd}{\tikzexternaldisable}
\AtEndEnvironment{tikzcd}{\tikzexternalenable}

\newtheorem{theorem}{Theorem}[section]
\newtheorem{lemma}[theorem]{Lemma}
\newtheorem{proposition}[theorem]{Proposition}
\newtheorem{corollary}[theorem]{Corollary}

\theoremstyle{definition}
\newtheorem{definition}[theorem]{Definition}
\newtheorem{construction}[theorem]{Construction}

\newtheorem{remark}[theorem]{Remark}
\newtheorem{example}[theorem]{Example}

\newtheorem{algorithm}[theorem]{Algorithm}

\DeclareMathOperator{\ZZ}{\mathbb{Z}}
\DeclareMathOperator{\QQ}{\mathbb{Q}}

\DeclareMathOperator{\KK}{\mathbb{K}}
\DeclareMathOperator{\TT}{\mathbb{T}}
\DeclareMathOperator{\PP}{\mathbb{P}}

\DeclareMathOperator{\N}{\mathbf{N}}
\DeclareMathOperator{\E}{\mathbf{E}}
\DeclareMathOperator{\F}{\mathbf{F}}
\DeclareMathOperator{\M}{\mathbf{M}}
\DeclareMathOperator{\K}{\mathbf{K}}
\renewcommand{\P}{\mathbf{P}}

\DeclareMathOperator{\Pic}{\mathrm{Pic}}
\DeclareMathOperator{\Cl}{\mathrm{Cl}}
\DeclareMathOperator{\im}{\mathrm{im}}

\DeclareMathOperator{\coker}{\mathrm{coker}}

\DeclareMathOperator{\cone}{\mathrm{cone}}
\DeclareMathOperator{\lin}{\mathrm{lin}}
\DeclareMathOperator{\rank}{\mathrm{rank}}

\DeclareMathOperator{\red}{\mathrm{red}}
\DeclareMathOperator{\fix}{\mathrm{fix}}

\DeclareMathOperator{\tors}{\mathrm{tors}}
\DeclareMathOperator{\sing}{\mathrm{sing}}
\DeclareMathOperator{\lcm}{\mathrm{lcm}}

\DeclareMathOperator{\iso}{\cong}

\newcommand{\vp}{\varphi}
\newcommand{\ol}{\overline}
\newcommand{\cal}{\mathcal}
\newcommand{\hatcal}[1]{\hat{\cal #1}}
\newcommand{\tildecal}[1]{\tilde{\cal #1}}
\newcommand{\barcal}[1]{\bar{\cal #1}}

\title{The Picard index of a surface with torus action}

\author{Justus Springer}

\address{Mathematisches Institut, Universit\"at T\"ubingen,
	Auf der Morgenstelle 10, 72076 T\"ubingen, Germany}

\email{justus.springer@uni-tuebingen.de}

\begin{document}

\begin{abstract}
	We consider normal rational projective surfaces with torus action
	and provide a formula for their Picard index, that means the
	index of the Picard group inside the divisor class group.
	As an application, we classify the log del Pezzo
	surfaces with torus action of Picard number one up to Picard index
	\( 10\,000 \).
\end{abstract}

\maketitle

\section{Introduction}

Consider a normal variety \( X \) defined over an algebraically
closed field \( \KK \) of characteristic zero. If~\( X \) is normal,
then the Picard group \( \Pic(X) \) embeds into the divisor class
group \( \Cl(X) \) as the subgroup consisting of the Cartier divisor
classes, and the \emph{Picard index} \( [\Cl(X):\Pic(X)] \) measures
the amount of non-invertible reflexive rank one sheaves on \( X \).
For rational normal projective surfaces admitting a (non-trivial)
action of the multiplicative group \( \KK^* \), we provide the
following formula, involving the torsion part \(
\Cl(X)^{\mathrm{tors}} \) and the local class groups \( \Cl(X,x) \),
hosting the Weil divisors modulo those being principal near \( x \in
X \).

\begin{theorem}
	\label{thm:picard_index_formula}
	The Picard index of a normal rational projective \( \KK^* \)-surface
	\( X \) is given by
	\[
		[\Cl(X) : \Pic(X)]\ =\
		\frac{1}{|\Cl(X)^{\mathrm{tors}}|} \prod_{x \in X} |\Cl(X,x)|.
	\]
\end{theorem}

Note that rationality forces our \( \KK^* \)-surface \( X \) to be \(
\QQ \)-factorial and \( \Cl(X) \) to be finitely generated, see for
instance~\cite{ArDeHaLa}*{Thm. 5.4.1.5}. Moreover, by normality of \(
X \), there are only finitely many singular points and these are the
only possible contributors of non-trivial local class groups. Thus,
all terms in our formula are indeed finite.

Beyond the \( \KK^* \)-surfaces, the formula trivially holds for all
smooth projective surfaces with a finitely generated and torsion free
divisor class group. As soon as we allow torsion, the r.h.s. is no
longer integral in the smooth case and thus the formula fails.
Concrete examples are the Enriques surfaces, having divisor class
group \( \ZZ^{10} \times \ZZ/2\ZZ \). A singular counterexample
without \( \KK^* \)-action is provided by the \( D_8 \)-singular log
del Pezzo surface of Picard number one: it is \( \QQ \)-factorial
with divisor class group \( \ZZ \times \ZZ / 2 \ZZ \) and doesn't
satisfy the formula, see Example~\ref{exm:d8_log_del_pezzo}.

Our motivation to consider the Picard index arises from the study of
log del Pezzo surfaces. Recall that these are normal projective
surfaces with an ample anticanonical divisor and at most finite
quotient singularities. The log del Pezzo surfaces form an infinite
class, which can be filtered into finite subclasses by further
conditions on the singularities. Common conditions are bounding the
Gorenstein index or the log terminality; for the state of the art we
refer to~\cites{AlNi, Nak, FuYa} in the general case and
to~\cites{HaHaHaSp, HaHaeSp, Hae} in the case of log del Pezzo
surfaces with \( \KK^* \)-action. The idea of filtering by the Picard
index has appeared in~\cite{HaHeSu}, where not-necessarily log
terminal Fano varieties with divisor class group \( \ZZ \) and a
torus action of complexity one have been considered. Here, we use
Theorem \ref{thm:picard_index_formula} to derive in Picard number one
suitable bounds on toric and non-toric log del Pezzo \( \KK^*
\)-surfaces and then present a classification algorithm. Explicit
results are obtained up to Picard index \( 1\,000\,000 \) in the
toric case and up to Picard index \( 10\,000 \) in the non-toric
case. The defining matrices for the surfaces are available under
\cite{ldp}, the toric ones up to Picard index \( 10\,000 \) and the
non-toric ones up to Picard index \( 2\,500 \).

\begin{theorem}
	\label{thm:classification_total}
	There are \( 1\,415\,486 \) families of log del Pezzo \( \KK^* \)-surfaces
	of Picard number one and Picard index at most \( 10\,000 \). Of those, \(
	68\,053 \) are toric and \( 1\,347\,433 \) are non-toric. The number of
	families for given Picard index develops as follows:

	\begin{center}
		\includegraphics{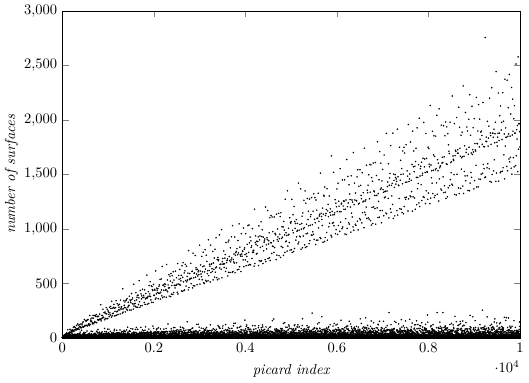}
	\end{center}

\end{theorem}

Let us give an idea of the proof of
Theorem~\ref{thm:picard_index_formula}. We use the approach
of~\cites{HaHe,HaSu}, encoding rational projective \( \mathbb{K}^* \)
-surfaces \( X \) via integral \( n \times r \) matrices~ \( P \) ,
see also \cite{ArDeHaLa}*{Sections 3.4, 5.4}. The columns \( v_1,
\ldots, v_r \) of \( P \) are the primitive ray generators of the fan
\( \Sigma \) of a toric variety~\( Z \) hosting our \( X \) as a
closed \( \mathbb{K}^* \)-invariant subvariety. As we will mention in
Proposition~\ref{prp:picard_index_eq_picard_index_toric_ambient}, a
basic feature of this encoding is that \( X \) inherits divisor class
group, Picard group and local class groups from its ambient toric
variety:
\[
	\Cl(X) \ = \ \Cl(Z), \qquad \Pic(X)\ =\ \Pic(Z), \qquad \Cl(X,x)\
	=\ \Cl(Z,x).
\]
This allows us to work entirely in the setting of toric varieties and
to provide in Proposition~\ref{prp:picard_index_of_toric_variety} a
first formula of the Picard index~\( \iota_{\mathrm{Pic}}(X) \),
which is in fact, independently from the embedded \( X \), valid for
any \( \QQ \)-factorial toric variety \( Z \):
\[
	\iota_{\mathrm{Pic}}(X)
	\ = \
	\iota_{\mathrm{Pic}}(Z)
	\ = \
	\frac{1}{\vert \coker(\hat P^*) \vert}
	\prod_{\sigma \in S} \vert \Cl(Z,z_\sigma) \vert.
\]
Here, \( S \subseteq \Sigma \) contains the cones \( \sigma \)
defining a closed torus orbit \( \TT^n \cdot z_{\sigma} \subseteq Z
\). The lattice homomorphism \( \hat P \) shows up when we describe
the Picard group inside the divisor class group as the intersection
over the kernels of the projections onto the local class groups of
the points \( z_\sigma \in Z \). In terms of the defining data, \(
\hat P \) is given as follows. Every cone \( \sigma \in \Sigma \)
defines lattices
\[
	{\textstyle
			N_{\sigma} \ := \ \lin_{\QQ}(\sigma) \cap \ZZ^n,
			\qquad \qquad
			F_{\sigma} \ := \ \underset{v_i\in\sigma}{\bigoplus} \ZZ\cdot f_{\sigma,i}.
		}
\]
Moreover, we have the homomorphisms \( P_{\sigma} \colon F_{\sigma}
\to N_{\sigma} \), sending the basis vector \( f_{\sigma,i} \in
F_\sigma \) to \( v_i \in N_\sigma \). These fit together to a
homomorphism
\[
	{\textstyle
			\underset{\sigma\in S}{\bigoplus}  P_{\sigma} \colon
			\underset{\sigma\in S}{\bigoplus} F_{\sigma}
			\ \to \
			\underset{\sigma\in S}{\bigoplus} N_{\sigma}.
		}
\]
The desired lattice homomorphism \(\hat P \colon \ker(\beta) \to
\ker(\alpha)\) is the (well defined) restriction of $\oplus
	P_{\sigma}$ to the kernels of
\[
	{\textstyle
			\beta \colon \underset{\sigma\in S}{\bigoplus} F_{\sigma} \to \ZZ^r,
			\quad
			f_{\sigma, i} \mapsto e_i,
			\qquad\qquad
			\alpha \colon
			\underset{\sigma\in S}{\bigoplus} N_{\sigma} \to \ZZ^n,
			\quad
			(v_{\sigma})_{\sigma\in S} \mapsto \underset{\sigma\in S}{\sum} v_{\sigma}.
		}
\]
At this point the serious and technical part of the proof of
Theorem~\ref{thm:picard_index_formula} begins, namely to verify the
identity
\[
	\vert \coker(\hat P^*) \vert
	\ = \
	\vert \Cl(Z)^{\mathrm{tors}} \vert.
\]
From toric geometry, we know \( \Cl(Z)^{\mathrm{tors}} = \coker(P)
\). In order to compare the cokernel orders, we give an explicit
matrix representation of \( \hat P \) and work with maximal minors on
both sides, see
Proposition~\ref{prp:gcd_maximal_minors_P_eq_gcd_maximal_minors_Phat}.
This part of the proof heavily depends on the specific shape of the
defining matrices \( P \) for \( \KK^* \)-surfaces. Once this is
done, Theorem~\ref{thm:picard_index_formula} follows from the
observation that the points \( x \in X \) with non-trivial local
class group correspond to points \( z_{\sigma} \) with \( \sigma
\in S \).

This article is structured as follows. In
Section~\ref{sec:picard_group_of_toric_variety}, we begin
the study of the
Picard index in a purely toric setting. In
Section~\ref{sec:cstar_surfaces}, we recall the combinatorial theory
of \( \KK^* \)-surfaces in terms of the matrices \( P \) mentioned
before. The technical part of the proof of
Theorem~\ref{thm:picard_index_formula} starts in
Section~\ref{sec:maximal_minors_P}, where we study maximal minors of
\( P \). In Section~\ref{sec:construction_of_phat}, we give an
explicit matrix representation of \( \hat P \) and in
Section~\ref{sec:maximal_minors_phat}, we study its maximal minors.
In Section~\ref{sec:proof_of_picard_index_formula}, we establish the
equality of the cokernel orders of~\( P \) and~\( \hat P^* \),
completing the proof of Theorem~\ref{thm:picard_index_formula}.
Finally, in Section~\ref{sec:log_del_pezzo_classification}, we
present the classification algorithm for log del Pezzo \( \KK^*
\)-surfaces of Picard number one with given Picard index, proving
Theorem~\ref{thm:classification_total}.

The author would like to thank Prof. J\"{u}rgen Hausen for his
valuable feedback and advice. Moreover, the author thanks the
referees for their very helpful comments and suggestions for
improving the presentation of the article.

\tableofcontents

\section{The Picard group of a toric variety}
\label{sec:picard_group_of_toric_variety}

In this section, we develop our approach to the Picard group of toric
varieties, which yields in
Corollary~\ref{cor:picard_group_free_of_alpha_surjective} a criterion
for torsion-freeness and in
Proposition~\ref{prp:picard_index_of_toric_variety} a first formula
involving the Picard index and local class groups. The reader is
assumed to be familiar with the basics of toric geometry
\cites{CoLiSc, Ful}.

\begin{construction}
	\label{cns:local_class_group_of_toric_variety}

	Let \( Z = Z_{\Sigma} \) be a toric variety coming from a fan \(
	\Sigma \) in the lattice \( N := \ZZ^r \). We assume \( \Sigma \) to
	be non-degenerate, i.e. its primitive ray generators \( v_1, \dots,
	v_n \) span \( \QQ^r \) as a convex cone. We allow \( \Sigma \) to be
	non-complete. With \( F := \ZZ^n \), we have the generator map
	\[
		P \colon F \to N, \qquad e_i \mapsto v_i.
	\]
	To any maximal cone \( \sigma = \cone(v_{i_1}, \dots,
	v_{i_{n_\sigma}}) \in \Sigma_{\max} \), we associate the lattices
	\[
		N_\sigma\ :=\ \lin_{\QQ}(\sigma) \cap N, \qquad \qquad F_\sigma\ :=\
		\ZZ^{n_\sigma}.
	\]
	We define the \emph{local generator map} associated to \( \sigma \)
	by
	\[
		P_\sigma \colon F_\sigma \to N_\sigma, \qquad e_j \mapsto v_{i_j}.
	\]
	With the inclusion \( \alpha_\sigma \colon N_\sigma \hookrightarrow N
	\) and the map \( \beta_\sigma \colon F_\sigma \to F \) sending \(
	e_j \) to \( e_{i_j} \), we obtain a commutative diagram
	\[
		\begin{tikzcd}
			F \ar[r, "P"]                                                      & N        \\
			F_\sigma \ar[u, hookrightarrow, "\beta_\sigma"] \ar[r, "P_\sigma"] & N_\sigma
			\ar[u, hookrightarrow, "\alpha_\sigma"].
		\end{tikzcd}
	\]
	Consider the dual lattices
	\[
		M := N^*, \qquad E := F^*,\qquad M_\sigma := N_\sigma^*,\qquad E_\sigma
		:= F_\sigma^*.
	\]
	Setting \( K := M/\im(P^*) \) and \( K_\sigma := M_\sigma /
	\im(P_\sigma^*) \), we obtain a map \( \pi_\sigma \colon K \to
	K_\sigma \) fitting into the commutative diagram with exact rows
	\[
		\begin{tikzcd}
			0 \ar[r]                 & M \ar[r, "P^*"] \ar[d,
			"\alpha_\sigma^*"]       & E \ar[r]
			\ar[d, "\beta_\sigma^*"] & K \ar[r] \ar[d, "\pi_\sigma"] & 0                                   \\
			0 \ar[r]                 & M_\sigma \ar[r, "P_\sigma^*"] & E_\sigma \ar[r] & K_\sigma \ar[r] &
			0.
		\end{tikzcd}
	\]
	By standard toric geometry, we have isomorphisms \( K \iso \Cl(Z) \)
	and \( K_\sigma \iso \Cl(U_\sigma) \), where \( U_\sigma \) is the
	affine toric chart associated to \( \sigma \). Moreover, the map \(
	\pi_\sigma \) corresponds to the restriction of divisor classes \(
	[D] \mapsto [D|_{U_\sigma}] \). In particular, its kernel consists of
	those divisor classes that are principal on \( U_\sigma \).

\end{construction}

\goodbreak

\begin{construction}
	\label{cns:picard_group_of_toric_variety}

	In the setting of Construction
	\ref{cns:local_class_group_of_toric_variety}, we define lattices \(
	\N \) and \( \F \) and a lattice homomorphism \( \P \colon \N \to \F
	\) by
	\[
		\N := \bigoplus_{\sigma \in \Sigma_{\max}} N_\sigma, \qquad
		\F := \bigoplus_{\sigma \in \Sigma_{\max}} F_\sigma, \qquad
		\P := \bigoplus_{\sigma \in \Sigma_{\max}} P_\sigma.
	\]
	Furthermore, we define lattice homomorphisms
	\begin{align*}
		\alpha \colon \N \to N, & \qquad N_\sigma \ni v \mapsto
		\alpha_\sigma(v),
		\\[3pt]
		\beta \colon \F \to F,  & \qquad F_\sigma \ni w \mapsto \beta_\sigma(w).
	\end{align*}
	Let \( \gamma \colon \hat N \to \N \) be a kernel of \( \alpha \) and \(
	\delta \colon \hat F \to \F \) be a kernel of \( \beta \). We obtain an
	induced map \( \hat P \colon \hat F \to \hat N \) making the following
	diagram commute:
	\[
		\begin{tikzcd}
			F \ar[r, "P"]                            & N                       \\
			\F \ar[r, "\P"] \ar[u, "\beta"]          & \N \ar[u, "\alpha"]     \\
			\hat F \ar[r, "\hat P"] \ar[u, "\delta"] & \hat N \ar[u, "\gamma"]
		\end{tikzcd}
	\]
	Now consider the dual lattices \( \M := \N^* \) and \( \E := \F^* \)
	as well as the abelian group \(\K := \bigoplus_{\sigma \in
		\Sigma_{\max}} K_\sigma \). We define the map
	\[
		\pi \colon K \to \K, \quad w \mapsto (\pi_\sigma(w))_{\sigma \in \Sigma_{\max}}.
	\]
	Setting \( \hat M := \M / \im(\alpha^*) \) and \( \hat E := \E /
	\im(\beta^*) \) as well as \( \hat K := \K / \im(\pi) \), we obtain a
	map \( \hat P' \colon \hat M \to \hat E \) fitting into the following
	commutative diagram with exact rows:
	\[
		\begin{tikzcd}
			0 \ar[r]      & M \ar[r, "P^*"] \ar[d, "\alpha^*"] & E \ar[r] \ar[d,
			"\beta^*"]    & K \ar[r]
			\ar[d, "\pi"] & 0                                                          \\
			0 \ar[r]      & \M \ar[r, "\P^*"]\ar[d]            & \E \ar[r] \ar[d] & \K
			\ar[r]\ar[d]  & 0                                                          \\
			              & \hat M \ar[r, "\hat P'"]           & \hat E \ar[r]
			              & \hat K \ar[r]                      & 0.
		\end{tikzcd}
	\]
\end{construction}

\begin{proposition}
	\label{prp:kernel_pi_eq_picard_group}
	In Construction \ref{cns:picard_group_of_toric_variety}, the map \( \beta \)
	is surjective and there is an exact sequence
	\[
		\begin{tikzcd}
			0 \ar[r] & \Pic(Z) \ar[r] & \hat M \ar[r, "\hat P'"] & \hat E \ar[r] & \hat K \ar[r] &
			0.
		\end{tikzcd}
	\]
	Moreover, if \( \alpha \) is surjective, \( \hat M \) is torsion-free
	and \( \hat P' = \hat P^* \).
\end{proposition}

\begin{proof}
	Every primitive generator of a ray of \( \Sigma \) is a generator of
	some maximal cone. This implies that \( \beta \) is surjective, hence
	\( \beta^* \) is injective. As a subgroup of \( K\), the Picard group
	\( \Pic(Z) \) consists of the Cartier divisor classes, i.e. those
	that are principal on all affine toric charts \( U_{\sigma} \) for \(
	\sigma \in \Sigma_{\max} \). This means
	\[
		\Pic(Z) = \bigcap_{\sigma \in \Sigma_{\max}} \ker(\pi_\sigma) =
		\ker(\pi).
	\]
	Applying the snake lemma to the lower diagram of Construction
	\ref{cns:picard_group_of_toric_variety}, gives the exact sequence of
	the Proposition. The last statement is clear.
\end{proof}

\begin{corollary}
	\label{cor:picard_group_free_of_alpha_surjective}

	Assume that in Construction \ref{cns:picard_group_of_toric_variety},
	the map \( \alpha \) is surjective. Then the Picard group \( \Pic(Z)
	\) is torsion-free.

\end{corollary}

\begin{remark}
	\label{rem:picard_group_free_of_maximal_dim_cone}
	Corollary \ref{cor:picard_group_free_of_alpha_surjective} generalizes the
	well-known fact that if \( Z \) has a toric fixed point,
	its Picard group is torsion-free. Indeed, having a toric fixed point means
	having a cone \( \sigma \in \Sigma \) of maximal dimension. This implies \( N_\sigma = N \), hence \( \alpha \)
	is surjective.
\end{remark}

\begin{definition}
	\label{def:picard_index}
	The \emph{Picard index} of a normal \( \QQ \)-factorial variety \( X \) is
	defined as
	\[
		\iota_{\Pic}(X)\ :=\ [\Cl(X) : \Pic(Z)].
	\]
\end{definition}

\begin{proposition}
	\label{prp:picard_index_of_toric_variety}

	Let \( Z = Z_{\Sigma} \) be a toric variety with a non-degenerate
	simplicial fan \( \Sigma \). In the notation of Construction
	\ref{cns:picard_group_of_toric_variety}, we have
	\[
		\iota_{\Pic}(Z)\ =\ \frac{1}{|\hat K|} \prod_{\sigma \in \Sigma_{\max}}
		|K_\sigma|.
	\]
\end{proposition}

\begin{proof}
	Recall that \( \Pic(Z) = \ker(\pi) \)
	and \( \hat K = \K / \im(\pi) \). Since \( \Sigma \) is simplicial, each
	\( K_\sigma \) is finite, hence so is \( \K \). We obtain
	\[
		\iota_{\Pic}(Z) = [K : \ker(\pi)] = |\im(\pi)| = \frac{|\K|}{|\hat K|} =
		\frac{1}{|\hat K|} \prod_{\sigma \in \Sigma_{\max}} |K_\sigma|.
	\]
\end{proof}

\section{Background on \texorpdfstring{\( \KK^* \)}{}-surfaces}
\label{sec:cstar_surfaces}

We recall the construction of rational $\KK^*$-surfaces $X$ in terms
of defining matrices~$P$ from~\cites{HaHe, HaSu, HaWr} and provide
the necessary facts around it and fix our notation for the subsequent
sections; see also~\cite{ArDeHaLa}*{Sections 3.4, 5.4}. As a basic
feature of this approach, $X$ comes embedded into a specific toric
variety $Z$ which determines a significant part of the geometry of
$X$. In particular and most relevant for us, $X \subseteq Z$ directly
inherits divisor class group, Picard group, local class groups and
hence also the Picard index; see
Proposition~\ref{prp:picard_index_eq_picard_index_toric_ambient}.
This will allow us to prove Theorem~\ref{thm:picard_index_formula}
entirely in terms of toric geometry.

We start by recalling some aspects of the geometry of \( \KK^*
\)-surfaces and their fixed points, the major part of which has been
developed in \cites{OrWa1, OrWa2, OrWa3}. A \emph{\( \KK^*
	\)-surface} is an irreducible, normal surface \( X \) coming with an
effective morphical action \( \KK^* \times X \to X \). Let \( X \) be
a projective \( \KK^* \)-surface. For each point \( x \in X \), the
orbit map \( t \mapsto t \cdot x \) extends to a morphism \( \vp_x
\colon \PP_1 \to X \). This allows one to define
\[
	x_0\ :=\ \vp_x(0), \qquad \qquad x_{\infty}\ :=\ \vp_x(\infty).
\]
The points \( x_0 \) and \( x_{\infty} \) are fixed points for the \(
\KK^* \)-action and they lie in the closure of the orbit \( \KK^*
\cdot x \). There are three types of fixed points: A fixed point is
called \emph{parabolic} (\emph{hyperbolic}, \emph{elliptic}), if it
lies in the closure of precisely one (precisely two, infinitely many)
non-trivial \( \KK^* \)-orbits. Hyperbolic and elliptic fixed points
are isolated, hence their number is finite. Parabolic fixed points
form a closed smooth curve with at most two connected components.
Every projective \( \KK^* \)-surface has a \emph{source} and a
\emph{sink}, i.e. two irreducible components of the fixed point set
\( F^+, F^- \subseteq X \) such that there exist non-empty \( \KK^*
\)-invariant open subsets \( U^+, U^- \subseteq X \) with
\[
	x_0 \in F^+ \text{ for all } x \in U^+, \qquad \qquad x_{\infty} \in F^- \text{
		for all } x \in U^-.
\]
The source either consists of a single elliptic fixed point or it is
a smooth curve of parabolic fixed points. The same holds for the
sink.

Before giving the general construction of \( \KK^* \)-surfaces in
terms of integral matrices, we start with an example, which will
reappear as a running example throughout the article.

\begin{example}
	\label{exm:running_example_toric_embedding}
	Consider the integral matrix
	\[
		P =
		\begin{bmatrix}
			v_{01} & v_{02} & v_{11} & v_{21}
		\end{bmatrix} =
		\begin{bmatrix}
			-1 & -1 & 8 & 0 \\
			-1 & -1 & 0 & 4 \\
			-1 & -2 & 7 & 3
		\end{bmatrix}.
	\]
	Let \( Z \) be the toric variety coming from the lattice fan \(
	\Sigma \) with maximal cones
	\[
		\sigma^+ := \cone(v_{01},v_{11},v_{21}), \quad \sigma^- :=
		\cone(v_{02},v_{11},v_{21}), \quad \tau_{01} := \cone(v_{01},v_{02}).
	\]
	We write \( U_1,U_2,U_3 \) for the coordinate functions of the acting
	torus \( \TT^3 \) of \( Z \) and set \( h := 1 + U_1 + U_2 \). Taking
	the closure of the zero locus defines an irreducible normal surface
	\[
		X := \overline{V(h)} \subseteq Z,
	\]
    with a \( \KK^* \)-action given on \( X \cap \TT^3 \) by \( t
    \cdot z = (z_1,z_2,tz_3) \).
	There are two elliptic fixed points \(
	x^{\pm} = z_{\sigma^{\pm}} \) and one hyperbolic fixed point \( \{
	x_{01} \} = X \cap \TT^3 \cdot z_{\tau_{01}} \). As we will see in
    Proposition~\ref{prp:picard_index_eq_picard_index_toric_ambient},
    divisor class group and Picard group can be computed via toric
    geometry:
	\[
		\Cl(X)\ =\ \ZZ \times \ZZ/4\ZZ, \qquad \Pic(X)\ =\
		\ZZ\cdot(15,\ol{1}) \subseteq \Cl(X).
	\]
	In particular, we have \( \iota_{\Pic}(X) = 60 \).

\end{example}

We come to the general construction of \( \KK^* \)-surfaces that will
be our working environment for the rest of this article. In a first
step, we produce our defining matrices \( P \).

\begin{construction}
	\label{cns:generator_matrix_cstar_surface}

	Fix positive integers \( r, n_0, \dots, n_r \). We start with
	integral vectors \(l_i = (l_{i1}, \dots, l_{in_i}) \in \ZZ_{\geq
		1}^{n_i} \) and \( d_i = (d_{i1}, \dots, d_{in_i}) \in \ZZ^{n_i}\)
	such that
	\[
		\gcd(l_{ij}, d_{ij}) = 1, \qquad
		\frac{d_{i1}}{l_{i1}} > \dots > \frac{d_{in_i}}{d_{in_i}}, \qquad
		i = 0, \dots, r.
	\]
	The building blocks for our defining matrices are
	\[
		L\ :=\
		\begin{bmatrix}
			-l_0   & l_1 & \dots  & 0      \\
			\vdots &     & \ddots & \vdots \\
			-l_0   & 0   & \dots  & l_r
		\end{bmatrix}, \qquad
		d\ :=\
		\begin{bmatrix}
			d_0 & d_1 & \dots & d_r
		\end{bmatrix}.
	\]
	According to the possible constellations of source and sink, we
	introduce four types of integral matrices:
	\[
		\begin{array}{crclcrcl}
			\text{(ee)}                  & P & = &
			\begin{bmatrix}
				L \\
				d
			\end{bmatrix}, \qquad \qquad &
			\text{(pe)}                  & P & = &
			\begin{bmatrix}
				L & 0 \\
				d & 1
			\end{bmatrix},                         \\\\
			\text{(ep)}                  & P & = &
			\begin{bmatrix}
				L & 0  \\
				d & -1
			\end{bmatrix}, \qquad \qquad &
			\text{(pp)}                  & P & = &
			\begin{bmatrix}
				L & 0 & 0  \\
				d & 1 & -1
			\end{bmatrix}.
		\end{array}
	\]
	With the canonical basis vectors \( e_1, \dots, e_{r+1} \) of \(
	\ZZ^{r+1} \) and \( e_0 := -(e_1 + \dots + e_r) \), the columns of \(
	P \) are
	\[
		v_{ij}\ :=\ l_{ij} e_i + d_{ij} e_{r+1}, \qquad v^+\ :=\ e_{r+1}, \qquad
		v^-\ :=\ -e_{r+1},
	\]
	where \( i = 0, \dots, r \) and \( j = 1, \dots, n_i \). We call \( P
	\) a \emph{defining matrix}, if its columns generate \( \QQ^{r+1} \)
	as a convex cone.
\end{construction}

Next, we will construct the fan for the ambient toric varieties for
our \( \KK^* \)-surfaces.
\begin{construction}
	\label{cns:ambient_toric_variety_cstar_surface}

	Let \( P \) be a defining matrix. Setting \( v_{i0} := v^+ \) and \(
	v_{in_i+1} := v^-\) for all \( i \), we define the cones
	\[
		\sigma^+ := \cone(v_{01}, \dots, v_{r1}), \qquad \sigma^- :=
		\cone(v_{0n_0}, \dots, v_{rn_r}),
	\]
	\[
		\tau_{ij} := \cone(v_{ij}, v_{ij+1}), \quad \text{for } i = 0, \dots, r
		\text{ and } j = 0, \dots, n_i.
	\]
	According to the type of \( P \), we define \( \Sigma \) to be the
	fan with the following maximal cones:
	\[
		\begin{array}{cccccc}
			\text{(ee)}                          & \{ \sigma^+ \}                    & \cup                   & \{ \tau_{i1}, \dots,
			\tau_{in_i-1}\ ;\ i = 0, \dots, r \} & \cup
			                                     & \{ \sigma^- \},                                                                   \\[3pt]
			\text{(pe)}                          & \{ \tau_{00}, \dots, \tau_{r0} \} & \cup                   & \{ \tau_{i1}, \dots,
			\tau_{in_i-1}\ ;\ i = 0, \dots, r \} & \cup
			                                     & \{ \sigma^- \},                                                                   \\[3pt]
			\text{(ep)}                          & \{ \sigma^+ \}                    & \cup                   & \{ \tau_{i1}, \dots,
			\tau_{in_i-1}\ ;\ i = 0, \dots, r \} & \cup                              & \{ \tau_{0n_0}, \dots,
			\tau_{rn_r} \},                                                                                                          \\[3pt]
			\text{(pp)}                          & \{ \tau_{00}, \dots,
			\tau_{r0} \}                         & \cup                              & \{ \tau_{i1}, \dots,
			\tau_{in_i-1}\ ;\ i = 0, \dots, r \} & \cup                              & \{ \tau_{0n_0}, \dots,
			\tau_{rn_r} \}.                                                                                                          \\[3pt]
		\end{array}
	\]
	Note that \( \Sigma \) is a non-degenerate simplicial lattice fan in
	\( \ZZ^{r+1} \). However, it is in general not complete.

\end{construction}

\begin{construction}
	\label{cns:cstar_surface_from_data}

	Let \( P \) be a defining matrix. Consider the toric variety \( Z =
	Z_{\Sigma} \), where \( \Sigma \) is as in Construction
	\ref{cns:ambient_toric_variety_cstar_surface}. Let \( U_1, \dots,
	U_{r+1} \) be the coordinate functions on the acting torus \(
	\TT^{r+1} \) of \( Z \). Fix pairwise different \( 1 = \lambda_2,
	\dots, \lambda_r \in \KK^* \) and set
	\[
		h_i := \lambda_i + U_1 + U_i, \qquad \qquad i = 2, \dots, r.
	\]
	Passing to the closure of the common set of zeroes of \( h_2, \dots,
	h_r \), we obtain an irreducible rational normal projective surface
	\[
		X(P) := \ol{V(h_2, \dots, h_r)} \subseteq Z.
	\]
	Since the \( h_i \) do not depend on the last coordinate \( U_{r+1}
	\) of \( \TT^{r+1} \), we get an effective \( \KK^* \)-action on \(
	X(P) \) as a subtorus of \( \TT^{r+1} \) by
	\[
		t \cdot x := (1, \dots, 1, t) \cdot x.
	\]
\end{construction}

\begin{example}
	\label{exm:example_as_instance_of_construction}
	The \( \KK^* \)-surface from
	Example~\ref{exm:running_example_toric_embedding} arises from
	Constructions~\ref{cns:generator_matrix_cstar_surface}
	to~\ref{cns:cstar_surface_from_data} by setting \( (n_0,n_1,n_2) =
	(2,1,1) \) as well as \( (l_0,l_1,l_2) = ((1,1),(8),(4)) \) and \(
	(d_0,d_1,d_2) = ((-1,-2), (7), (3)) \) and choosing the type (ee).
\end{example}

\begin{remark}
	\label{rem:fixed_points_cstar_surface_of_data}

	Consider a \( \KK^* \)-surface \( X = X(P) \). Let \( Z = Z_{\Sigma}
	\) be the ambient toric variety with acting torus \( T = \TT^{r+1}
	\). The fixed points of \( X \) are given as follows. For every \(
	\tau_{ij} \in \Sigma_{\max} \), the associated toric orbit \( T \cdot
	z_{\tau_{ij}} \) intersects \( X \) in a fixed point
	\[
		\{x_{ij}\}\ =\ X\ \cap\ T \cdot z_{\tau_{ij}}.
	\]
	If \( 1 \leq j \leq n_i-1 \), the fixed point \( x_{ij} \) is
	hyperbolic and all hyperbolic fixed points arise this way. For \( j
	\in \{0, n_i\} \), the fixed point \( x_{ij} \) is parabolic.
	According to the type of \( P \), we have the following.
	\begin{itemize}
		\item[(ee)] There are two elliptic fixed points \( x^+ = z_{\sigma^+} \)
			and \( x^- = z_{\sigma^-} \) and no parabolic fixed points.
		\item[(pe)] There is one elliptic fixed point \( x^- = z_{\sigma^-} \).
			There are parabolic fixed points \( x_{i0} \in F^+ \) and all
			parabolic fixed points in \( F^+ \backslash \{x_{00}, \dots, x_{r0}\} \)
			are smooth.
		\item[(ep)] There is one elliptic fixed point \( x^+ = z_{\sigma^+} \).
			There are parabolic fixed points \( x_{in_i} \in F^- \) and all
			parabolic fixed points in \( F^- \backslash
			\{x_{0n_0}, \dots, x_{rn_r}\} \) are smooth.
		\item[(pp)] There are no elliptic fixed points. There are parabolic
			fixed points \( x_{i0} \in F^+ \) and \( x_{in_i} \in F^- \) and all
			parabolic fixed points in \( F^+
			\backslash \{x_{00}, \dots, x_{r0}\} \) and \( F^- \backslash
			\{x_{0n_0}, \dots, x_{rn_r} \} \) are smooth.
	\end{itemize}
\end{remark}

\begin{theorem}[See {\cite{ArDeHaLa}*{Thm. 5.4.1.5}}]
	\label{thm:cstar_surface_from_data_completeness}
	Every rational projective \( \KK^* \)-surface is isomorphic to a \( \KK^*
	\)-surface \( X(P) \) arising from Construction
	\ref{cns:cstar_surface_from_data}.
\end{theorem}

\begin{proposition}
	\label{prp:picard_index_eq_picard_index_toric_ambient}

	Let \( X = X(P) \subseteq Z \) arise from
    Construction~\ref{cns:cstar_surface_from_data}. Then we have
	\begin{enumerate}[\normalfont(i)]
		\item \( \Cl(X) \iso \Cl(Z) \),
		\item \( \Pic(X) \iso \Pic(Z) \),
		\item \( \Cl(X, x) \iso \Cl(Z,x) \) for all \( x \in X \),
		\item \( \iota_{\Pic}(X) = \iota_{\Pic}(Z) \).
	\end{enumerate}
\end{proposition}
\begin{proof}
    By~\cite{ArDeHaLa}*{Prop. 5.4.1.8}, the embedding \( X \subseteq Z
    \) is the canonical toric embedding in the sense
    of~\cite{ArDeHaLa}*{Sec. 3.2.5}. Now use~\cite{ArDeHaLa}*{3.2.5.4,
    3.3.1.7}.
\end{proof}

\section{Maximal minors of \texorpdfstring{\( P \)}{P}}
\label{sec:maximal_minors_P}

As a first step towards the equality of the cokernel orders of $P$
and $\hat P^*$, we study in this section the set $M(P)$ of maximal
minors of our defining matrices $P$. We figure out relations among
the maximal minors of $P$ and, based on that, show that $\gcd(M(P))$
equals $\gcd(M'(P))$ for a proper, more accessible subset $M'(P)
	\subset M(P)$; see
Proposition~\ref{prp:gcd_maximal_minors_P_eq_gcd_maximal_minors_prime}.
First, we introduce a suitable lattice basis for~$\mathbb{Z}^n$,
reflecting the block structure of $P$.

\begin{construction}
	\label{cns:generator_matrix_cstar_surface_old}

	Let \( P \) be a defining matrix. Set \( n := n_0 + \dots + n_r \).
	Then \( P \) is an integral \( (r+1) \times (n+m) \)-matrix, where \(
	m \in \{0,1,2\} \). Write \( e_1, \dots, e_{r+1} \) for the canonical
	basis vectors of \( \ZZ^{r+1} \) and set \( u := e_{r+1} \). Then we
	define
	\[
		N \ := \ \ZZ^{r+1} \ = \ \ZZ e_1 \oplus \dots \oplus \ZZ e_r \oplus \ZZ
		u.
	\]
	For each \( i = 0, \dots, r \), we introduce a set of symbols to be
	used as a lattice basis:
	\[
		\begin{array}{lclclc>{\displaystyle}l}
			\cal F_i                                         & := & \{f_{i1}, \dots, f_{in_i}\},                    & \qquad \qquad & F_i & := & \ZZ f_{i1}
			\oplus \dots \oplus \ZZ f_{in_i}\ \iso \ \ZZ^{n_i},                                                                                             \\[8pt]
			\cal F'                                          & := &
			\begin{cases}
				\emptyset,    & \text{(ee)} \\
				\{f^+\},      & \text{(pe)} \\
				\{f^-\},      & \text{(ep)} \\
				\{f^+, f^-\}, & \text{(pp)}
			\end{cases}, &    &
			F'                                               & := & \bigoplus_{f \in \cal F'} \ZZ f \ \iso \ \ZZ^m.
		\end{array}
	\]
	Setting
	\[
		\cal F := \cal F_0 \cup \dots \cup \cal F_r \cup \cal F', \qquad \qquad
		F \ := \ F_0 \oplus \dots \oplus F_r \oplus F',
	\]
	we obtain an isomorphism \( F \iso \ZZ^{n+m} \). Set \( e_0 := -(e_1
	+ \dots + e_r) \). Then we can view \( P \) as a lattice map given by
	\[
		P \colon F \to N, \quad
		\begin{cases}
			f_{ij} \mapsto v_{ij} := l_{ij} e_i + d_{ij} u &                             \\
			f^+ \mapsto u,                                 & \text{ if } f^+ \in \cal F' \\
			f^- \mapsto -u,                                & \text{ if } f^- \in
			\cal F'.
		\end{cases}
	\]

\end{construction}

\begin{definition}
	\label{def:maximal_minors_P_cstar_surface}
	Let a subset \( A \subseteq \cal F \) with \( |A| = r+1 \) be given. Then we
	have a sublattice \( F_A := \bigoplus_{f \in A} \ZZ \cdot f \subseteq F \) and
	an induced map \( P_A \colon F_A \to N \) as in the commutative diagram
	\[
		\begin{tikzcd}
			F \ar[r, "P"] & N \\
			F_A \ar[u, hookrightarrow] \ar[ur, "P_A" {below, yshift=-0.5ex}].
		\end{tikzcd}
	\]
	We call \( |\det(P_A)| \in \ZZ_{\geq 0} \) the \emph{maximal minor of
		\( P \) associated to \( A \)}. The set of maximal minors of \( P \)
	is
	\[
		M(P)\ :=\ \{ |\det(P_A)|\ ;\ A \subseteq \cal F,\ |A| = r+1 \}.
	\]
\end{definition}

\begin{example}
	\label{exm:running_example_maximal_minors}
	With \( P \) as in
	Example~\ref{exm:running_example_toric_embedding}, we have \( \cal
	F = \{f_{01}, f_{02}, f_{11}, f_{21} \} \). Setting \( A_{ij} :=
	\cal F \backslash \{f_{ij}\} \), the maximal minors of \( P
	\) are
	\begin{align*}
		\vert \det(P_{A_{01}}) \vert & = \vert d_{02}l_{11}l_{21} +
		l_{02}d_{11}l_{21}+l_{02}l_{11}d_{21} \vert = 12,                  \\
		\vert \det(P_{A_{02}}) \vert & = \vert d_{01}l_{11}l_{21} +
		l_{01}d_{11}l_{21}+l_{01}l_{11}d_{21}\vert = 20,                   \\
		\vert \det(P_{A_{11}}) \vert & = \vert
        l_{21}(d_{01}l_{02}-l_{01}d_{02}) \vert = 4, \\
		\vert \det(P_{A_{21}}) \vert & = \vert
        l_{11}(d_{01}l_{02}-l_{01}d_{02}) \vert = 8.
	\end{align*}
	Hence, we have \( M(P) = \{12,20,4,8\} \).

\end{example}

The following lemma gives a vanishing criterion for maximal minors of
\( P \). In particular, it will allow us to describe the nonzero
maximal minors of \( P \) explicitly.

\begin{lemma}
	\label{lem:maximal_minor_P_vanishing_cstar_surface}
	Let \( A \subseteq \mathcal{F} \) with \( |A| = r + 1 \). Assume that there
	exist \( 0 \leq i_0 < i_1 \leq r \) such that \( A \cap \mathcal{F}_{i_0} = A
	\cap \mathcal{F}_{i_1} = \emptyset \). Then \( \det(P_A) = 0 \).
\end{lemma}

\begin{proof}
	Consider the dual bases \( \{f^*_{ij}\} \) and \( \{e_1^*, \dots, e_r^*,
	u^*\} \) of \( F^* \) and \( N^* \) respectively. Then we have
	\[
		P^*(e_i^*) = \left( \sum_{j=1}^{n_i} l_{ij}
		f^*_{ij} \right) - \left( \sum_{j=1}^{n_0} l_{0j} f^*_{0j} \right).
	\]
	If \( i_0 = 0 \), this implies \( P^*_A(e_{i_1}^*) = 0 \). If \( i_0
	\) and \( i_1 \) are both nonzero, we have
	\[
		P^*_A(e_{i_0}^*) = - \left( \sum_{j=1}^{n_0} l_{0j} f^*_{0j} \right) =
		P^*_A(e_{i_1}^*).
	\]
	Thus in both cases, \( \det(P_A) = \det(P_A^*) = 0 \).
\end{proof}

\goodbreak

\begin{definition}
	\label{def:maximal_minors_mu_nu}
	\begin{enumerate}[\normalfont(i)]
		\item Let numbers \( 1 \leq j_i \leq n_i \) for all \( i = 0, \dots, r \)
		      be given. We define
		      \[
			      \mu(j_0, \dots, j_r)\ :=\ \sum_{i_0=0}^{r} d_{i_0 j_{i_0}} \prod_{i
				      \neq i_0} l_{ij_i}, \qquad \qquad \hat \mu\ :=\ \mu(n_0, \dots, n_r).
		      \]
		\item Let \( 0 \leq i \leq r \) and \( 0 \leq j, j' \leq n_i \) be given.
		      We define
		      \[
			      \nu(i, j, j')\ :=\ l_{ij} d_{ij'} - l_{ij'} d_{ij}, \qquad \qquad
			      \hat \nu(i, j)\ :=\ \nu(i, j, n_i).
		      \]
	\end{enumerate}
\end{definition}

\begin{lemma}
	\label{lem:nonzero_maximal_minors_P_of_intersection_Fp_empty}
	Let \( A \subseteq \mathcal{F} \) with \( |A| = r + 1 \) such that \( \det(P_A)
	\neq 0 \). Assume that \( A \cap \cal F' = \emptyset \) holds.
	Then either (i) or (ii) must hold.
	\begin{enumerate}[\normalfont(i)]
		\item We have \( |A \cap \mathcal{F}_i| = 1 \) for all \( i = 0, \dots, r
		      \) and \(|\det(P_A)| = |\mu(j_0, \dots, j_r)| \) for some numbers \(
		      1 \leq j_i \leq n_i \).
		\item There exist \( 0 \leq i_0, i_1 \leq r \) with
		      \[
			      |A \cap \mathcal{F}_i| =
			      \begin{cases}
				      2, & i = i_0          \\
				      0, & i = i_1          \\
				      1, & \text{otherwise}
			      \end{cases}
		      \]
		      and we have \(|\det(P_A)| = |\nu(i_0, j_{i_0}, j_{i_0}')| \prod_{i
			      \neq i_0,i_1} l_{ij_i}\) for some numbers \( 1 \leq j_i \leq n_i \)
		      as well as \( 1 \leq j_{i_0}' \leq n_{i_0} \).
	\end{enumerate}
\end{lemma}

\begin{proof}
	Lemma \ref{lem:maximal_minor_P_vanishing_cstar_surface} implies that there
	is at most one \( i = 0, \dots, r \) with \( A \cap \cal F_i = \emptyset \).
	Since \( A \cap \cal F' = \emptyset \) and \( |A| = r+1 \), the
	conditions on \( |A \cap \cal F_i| \) follow. The expressions for \( \det(P_A)
	\) then come from cofactor expansion.
\end{proof}

\begin{example}
	\label{exm:running_example_minor_structure}
	Applying
	Lemma~\ref{lem:nonzero_maximal_minors_P_of_intersection_Fp_empty}
	to Example~\ref{exm:running_example_maximal_minors}, we see that
	the minors associated to \( A_{01} \) and \( A_{02} \) satisfy
	(i), while those associated to \( A_{11} \) and \( A_{21} \)
	satisfy~(ii).
\end{example}

\begin{lemma}
	\label{lem:nonzero_maximal_minors_P_of_intersection_Fp_nonempty}
	Let \( A \subseteq \cal F \) with \( |A| = r+1 \) such that \( \det(P_A)
	\neq 0 \). Assume that \( A \cap \cal F' \neq \emptyset \) holds.
	Then there exists an \( i_1 = 0, \dots, r \) such that
	\[
		|A \cap \cal F_i| =
		\begin{cases}
			0, & i = i_1          \\
			1, & \text{otherwise}
		\end{cases}
	\]
	and we have \( |\det(P_A)| = \prod_{i \neq i_1} l_{ij_i} \).
\end{lemma}

\begin{proof}
	Note that \( |A \cap \cal F'| = 2 \) cannot hold, since this would
	imply \( \det(P_A) = 0 \). Hence \( |A \cap \cal F'| = 1 \). Lemma
	\ref{lem:maximal_minor_P_vanishing_cstar_surface} forces the condition on \(
	|A \cap \cal F_i| \). Cofactor expansion gives the expression for \(
	\det(P_A) \).
\end{proof}

\begin{lemma}
	\label{lem:nonzero_maximal_minors_mu_nu_relations}
	\begin{enumerate}[\normalfont(i)]
		\item Let \( i = 0, \dots, r \) and \( 1 \leq j,j' \leq n_i \). There exist
		      integers \( \alpha, \beta \in \ZZ \) such that
		      \[
			      \nu(i,j,j')\ =\ \alpha \hat \nu(i,j) - \beta \hat \nu(i,j').
		      \]
		\item Let \( i = 0, \dots, r \) and \( j_i = 1, \dots, n_i \) for all \( i
		      \). There exist integers \( \beta_i \) such that
		      \begin{align*}
			      \mu(j_0, \dots, j_r) \ = \  &
			      \beta_{i_0}\mu(j_0, \dots, j_{i_0-1}, n_{i_0}, j_{i_0+1}, \dots,
			      j_r)                                                            \\
			                                  & + \sum_{i_1 \neq i_0} \beta_{i_1}
			      \hat \nu(i_0, j_{i_0}) \prod_{i \notin \{i_0, i_1\}} l_{i
					      j_i}.
		      \end{align*}
	\end{enumerate}
\end{lemma}

\begin{proof}
	We show (i). Since \( \gcd(l_{in_i}, d_{in_i}) = 1 \), we find integers
	\( x, y \in \ZZ \) such that \( x l_{in_i} + y d_{in_i} = 1 \). Set \(
	\alpha := x l_{ij} + y d_{ij} \) and \( \beta := x l_{ij'} + y d_{ij'} \).
	Then the rows of \( 3 \times 3 \)-matrix
	\[
		\begin{bmatrix}
			l_{ij} & l_{ij'} & l_{in_i} \\
			d_{ij} & d_{ij'} & d_{in_i} \\
			\alpha & \beta   & 1
		\end{bmatrix}
	\]
	are linearly dependent, hence its determinant vanishes. Cofactor
	expansion by the last row gives the desired equality. For (ii), we
	consider exemplarily the case \( i_0 = 0 \). Since \( \gcd(l_{0n_0},
	d_{0n_0}) = 1 \), we find \( x, y \in \ZZ \) such that \( -x l_{0n_0}
	+ y d_{0n_0} = 1 \). Now set
	\[
		\beta_0 := -x l_{0 j_0} + y d_{0 j_0}, \qquad \qquad \beta_1 := x l_{1
				j_1} + y d_{1 j_1},
	\]
	as well as \( \beta_i := yd_{ij_i} \) for \( i = 2, \dots, r \). Then
	the rows of the \( (r+2) \times (r+2) \) matrix
	\[
		\begin{bmatrix}
			-l_{0j_0} & -l_{0n_0} & l_{1 j_1} & 0         & \dots  &
			0
			\\
			-l_{0j_0} & -l_{0n_0} & 0         & l_{2 j_2} & \dots  &
			0
			\\
			\vdots    & \vdots    & \vdots    &           & \ddots &
			\vdots
			\\
			-l_{0j_0} & -l_{0n_0} & 0         & \dots     & \dots  &
			l_{r j_r}
			\\
			d_{0j_0}  & d_{0n_0}  & d_{1j_1}  & d_{2j_2}  & \dots  &
			d_{r j_r}
			\\
			\beta_0   & 1         & \beta_1   & \beta_2   & \dots  &
			\beta_r
		\end{bmatrix}
	\]
	are linearly dependent, hence its determinant vanishes. Cofactor
	expansion by the last row and adjusting the signs of the \( \beta_i
	\) gives the desired equality.
\end{proof}

\begin{definition}
	\label{def:maximal_minors_prime_cstar_surface_ee}
	According to the type of \( P \), we define the set
	\[
		M'(P) :=
		\begin{cases}
			\displaystyle
			\{ |\hat \mu| \} \cup \left\{ |\hat \nu(i_0, j_{i_0})| \prod_{i \neq
				i_0, i_1} l_{i j_i}\ ;\
			\begin{aligned}
				 & 0 \leq i_0, i_1 \leq r,\text{ with } i_0 \neq i_1, \\
				 & 1 \leq j_i \leq n_i \text{ for all } i \neq i_1
			\end{aligned} \right\}, \quad
			\text{(ee)} \\
			\displaystyle
			\left\{ \prod_{i \neq i_1} l_{ij_i}\ ;\
			0 \leq i_1 \leq r,\ 1 \leq j_i \leq n_i \text{ for all } i \neq i_1
			\right\}, \quad
			\text{(pe), (ep), (pp).}
		\end{cases}
	\]
\end{definition}

\begin{proposition}
	\label{prp:gcd_maximal_minors_P_eq_gcd_maximal_minors_prime}
	We have \( \gcd(M(P)) = \gcd(M'(P)) \).
\end{proposition}

\begin{proof}
	Consider first the case (ee). Then \( \cal F' = \emptyset \), hence
	Lemma \ref{lem:nonzero_maximal_minors_P_of_intersection_Fp_empty}
	implies \( M'(P) \subseteq M(P) \). This shows that \( \gcd(M(P)) \)
	divides \( \gcd(M'(P)) \). On the other hand, by repeated application
	of Lemma \ref{lem:nonzero_maximal_minors_mu_nu_relations}, we see
	that every maximal minor of \( P \) can be written as a \( \ZZ
	\)-linear combination of elements of \( M'(P) \). This implies that
	\( \gcd(M'(P)) \) divides \( \gcd(M(P)) \). For cases (pe), (ep) and
	(pp), Lemma
	\ref{lem:nonzero_maximal_minors_P_of_intersection_Fp_nonempty}
	implies \( M'(P) \subseteq M(P) \), hence \( \gcd(M(P)) \) divides \(
	\gcd(M'(P)) \). On the other hand, Lemma
	\ref{lem:nonzero_maximal_minors_P_of_intersection_Fp_empty} implies
	that every maximal minor of \( P \) is a \( \ZZ \)-linear combination
	of elements of \( M'(P) \), hence \( \gcd(M'(P)) \) divides \(
	\gcd(M(P)) \).
\end{proof}

\begin{example}
	\label{exm:running_example_maximal_minors_prime}
	In Example~\ref{exm:running_example_maximal_minors}, we can apply
	Lemma~\ref{lem:nonzero_maximal_minors_mu_nu_relations} (ii) to
	express \( |\det(P_{A_{02}})| \) as a \( \ZZ
	\)-linear combination of the other maximal minors. Hence we have
	\( M'(P) = M(P) \backslash \{|\det(P_{A_{02}})|\} = \{12,8,4\} \).
\end{example}

\section{The construction of \texorpdfstring{\( \hat P \)}{P}}
\label{sec:construction_of_phat}

The aim of this section is to provide an explicit matrix
representation of the map~\( \hat P \) from
Construction~\ref{cns:picard_group_of_toric_variety}, which will be
used to describe its maximal minors in
Section~\ref{sec:maximal_minors_phat}.

\begin{construction}
	\label{cns:alpha_beta_cstar_surface}

	Let \( \Sigma \) be the fan of an ambient toric variety of a \( \KK^*
	\)-surface, as defined in Construction
	\ref{cns:ambient_toric_variety_cstar_surface}. Consider the lattices
	\( F_{\sigma} := \ZZ^{n_{\sigma}} \) and \( N_{\sigma} :=
	\lin_{N_{\QQ}}(\sigma) \cap N \) from Construction
	\ref{cns:local_class_group_of_toric_variety}. We have \( n_{\sigma^+}
	= n_{\sigma^-} = r+1 \) and \( n_{\tau_{ij}} = 2 \). Moreover, we
	have \( N_{\sigma^+} = N_{\sigma^-} = N = \ZZ^{r+1} \) and \(
	N_{\tau_{ij}} = \ZZ e_i + \ZZ u \subseteq N \). We will work with the
	identifications
	\[
		\begin{array}{lclclcl}
			F_{\sigma^+}  & \iso & \ZZ f_{01}^+ \oplus \dots \oplus \ZZ f_{r1}^+, &
			\qquad \qquad &
			N_{\sigma^+}  & \iso & \ZZ e_1^+ \oplus \dots \oplus \ZZ e_r^+ \oplus
			\ZZ u^+,                                                                       \\[5pt]
			F_{\tau_{ij}} & \iso & \ZZ f_{ij}^- \oplus \ZZ f_{ij+1}^+,            & \qquad
			\qquad        &
			N_{\tau_{ij}} & \iso & \ZZ e_{ij} \oplus \ZZ u_{ij},                           \\[5pt]
			F_{\sigma^-}  & \iso & \ZZ f_{0n_0}^- \oplus \dots \oplus \ZZ
			f_{rn_r}^-,   &
			\qquad \qquad &
			N_{\sigma^-}  & \iso & \ZZ e_1^- \oplus \dots \oplus \ZZ e_r^- \oplus
			\ZZ u^-.                                                                       \\[5pt]
		\end{array}
	\]
	Then according to the type of \( P \), a lattice basis of \( \N \) is
	given by
	\[
		\begin{array}{cccccc}
			\text{(ee)} & \{e_1^+, \dots, e_r^+, u^+\}                  & \cup & S & \cup & \{ e_1^-,
			\dots, e_r^-, u^- \},                                                                     \\[5pt]
			\text{(pe)} & \{ e_{i0}, u_{i0}\ ;\ i = 0, \dots, r \}      & \cup & S & \cup & \{ e_1^-,
			\dots, e_r^-, u^- \},                                                                     \\[5pt]
			\text{(ep)} & \{e_1^+, \dots, e_r^+, u^+\}                  & \cup & S &
			\cup        & \{ e_{in_i}, u_{in_i}\ ;\ i = 0, \dots, r \},                               \\[5pt]
			\text{(pp)} & \{ e_{i0}, u_{i0}\ ;\ i = 0, \dots, r \}      & \cup & S &
			\cup        & \{ e_{in_i}, u_{in_i}\ ;\ i = 0, \dots, r \},                               \\
		\end{array}
	\]
	where \( S := \{ e_{ij}, u_{ij}\ ;\ i = 0, \dots, r,\ j = 1, \dots,
	n_i-1 \} \). A lattice basis of \( \F \) is given by
	\[
		\begin{array}{cccc}
			\text{(ee)} & \{ f_{ij}^-\ ;\ i = 0, \dots, r,\ j = 1, \dots, n_i \}
			            & \cup
			            & \{ f_{ij}^+\ ;\ i = 0, \dots, r,\ j = 1, \dots, n_i
			\},                                                                  \\[5pt]
			\text{(pe)} & \{ f_{ij}^-\ ;\ i = 0, \dots, r,\ j = 0, \dots, n_i \}
			            & \cup
			            & \{ f_{ij}^+\ ;\ i = 0, \dots, r,\ j = 1, \dots, n_i
			\},                                                                  \\[5pt]
			\text{(ep)} & \{ f_{ij}^-\ ;\ i = 0, \dots, r,\ j = 1, \dots, n_i \}
			            & \cup
			            & \{ f_{ij}^+\ ;\ i = 0, \dots, r,\ j = 1, \dots, n_i+1
			\},                                                                  \\[5pt]
			\text{(pp)} & \{ f_{ij}^-\ ;\ i = 0, \dots, r,\ j = 0, \dots, n_i \}
			            & \cup
			            & \{ f_{ij}^+\ ;\ i = 0, \dots, r,\ j = 1, \dots, n_i+1
			\}.
		\end{array}
	\]
	In particular, we have \( \rank(\F) = \rank(\N) = 2n+m(r+1) \). With
	respect to these bases, the maps \( \alpha \) and \( \beta \) from
	Construction~\ref{cns:picard_group_of_toric_variety} are
	\[
		\begin{array}{ll}
			\alpha \colon \N \to N, & e_{ij}, e_i^+, e_i^- \mapsto e_i, \qquad
			u_{ij}, u^+, u^- \mapsto u,
			\\[5pt]
			\beta \colon \F \to F,  & f_{ij}^+, f_{ij}^- \mapsto f_{ij}.
		\end{array}
	\]
	Setting \( e^{+}_0 := -(e^{+}_1 + \dots + e^{+}_r) \) and \( e^{-}_0
	:= -(e^{-}_1 + \dots + e^{-}_r) \), the maps \( P_{\sigma^{+}},
	P_{\sigma^-} \) and \( P_{\tau_{ij}} \) are then given as
	\[
		\begin{array}{ll}
			P_{\sigma^+} \colon F_{\sigma^+} \to N_{\sigma^-},    & f^+_{i1} \mapsto l_{i1}
			e^+_i + d_{i1} u^+,
			\\[5pt]
			P_{\sigma^-} \colon F_{\sigma^-} \to N_{\sigma^-},    &
			f^-_{in_i} \mapsto l_{in_i} e^-_i + d_{in_i} u^-,
			\\[5pt]
			P_{\tau_{ij}} \colon F_{\tau_{ij}} \to N_{\tau_{ij}}, &
			f^-_{ij} \mapsto l_{ij} e_{ij} + d_{ij} u_{ij}
			\\[5pt]
			                                                      &
			f^+_{ij+1} \mapsto l_{ij+1} e_{ij} + d_{ij+1} u_{ij}.
		\end{array}
	\]

\end{construction}

\begin{remark}
	\label{rem:picard_group_free_of_cstar_surface}

	Clearly, the map \( \alpha \) in
	Construction~\ref{cns:alpha_beta_cstar_surface} is surjective in all
	cases. Hence
	Corollary~\ref{cor:picard_group_free_of_alpha_surjective} implies
	that the Picard group of a projective rational \( \KK^* \)-surface is
	torsion-free.

\end{remark}

\begin{example}
	\label{exm:running_example_alpha_beta}

	With \( P \) as in Example~\ref{exm:running_example_toric_embedding},
	we have \( \{e_1^+,e_2^+,u^+,e_{01},u_{01},e_1^-,e_2^-,u^-\} \) as a
	basis for \( \N \) and \( \{f_{01}^+, f_{11}^+, f_{21}^+, f_{01}^-,
	f_{02}^+, f_{02}^-, f_{11}^-, f_{21}^- \} \) as a basis for \( \F \).
	The matrix representations of \( \alpha \) and \( \beta \) are
	\[
		\alpha\ =\
		\begin{bmatrix}
			1 & 0 & 0 & -1 & 0 & 1 & 0 & 0 \\
			0 & 1 & 0 & -1 & 0 & 0 & 1 & 0 \\
			0 & 0 & 1 & 0  & 1 & 0 & 0 & 1
		\end{bmatrix}, \qquad
		\beta\ =\
		\begin{bmatrix}
			1 & 0 & 0 & 1 & 0 & 0 & 0 & 0 \\
			0 & 0 & 0 & 0 & 1 & 1 & 0 & 0 \\
			0 & 1 & 0 & 0 & 0 & 0 & 1 & 0 \\
			0 & 0 & 1 & 0 & 0 & 0 & 0 & 1
		\end{bmatrix}.
	\]

\end{example}

\begin{construction}
	\label{cns:phat_cstar_surface}
	We continue in the setting of
	Construction~\ref{cns:alpha_beta_cstar_surface}. We will give
	explicit descriptions of the maps \( \gamma, \delta \) and \( \hat P \)
	from Construction \ref{cns:picard_group_of_toric_variety}. According to the
	type of \( P \), let us set for all \( i = 0, \dots, r \)
	\[
		n_i' :=
		\begin{cases}
			n_i - 1, & \text{(ee)},             \\
			n_i,     & \text{(pe), (ep), (pp)},
		\end{cases}
	\]
	\[
		\hatcal N_i := \{ \hat e_{ij}, \hat u_{ij}\ ;\ j = 1, \dots, n_i' \}, \qquad
		\qquad \hatcal F_i := \{\hat f_{ij}\ ;\ j = 1, \dots n_i \}.
	\]
	Now define the sets of symbols
	\[
		\begin{array}{lc>{\displaystyle}lclc>{\displaystyle}l}
			\tildecal N                                                                      & :=     &
			\begin{cases}
				\{\tilde e_1, \dots, \tilde e_r, \tilde u\}, & \text{(ee)} \\
				\emptyset,                                   & \text{(pe)} \\
				\emptyset,                                   & \text{(ep)} \\
				\{\tilde u_1, \dots, \tilde u_r, \tilde e\}, & \text{(pp)} \\
			\end{cases}, & \qquad &
			\barcal F                                                                        & :=     &
			\begin{cases}
				\emptyset,                                       & \text{(ee)} \\
				\{\hat f_i^-\ ;\ i = 1, \dots, r \},             & \text{(pe)} \\
				\{\hat f_i^+\ ;\ i = 1, \dots, r \},             & \text{(ep)} \\
				\{\hat f_i^+, \hat f_i^-\ ;\ i = 1, \dots, r \}, & \text{(pp)} \\
			\end{cases},                                                                                    \\\\
			\hatcal N                                                                        & :=     & \hatcal N_0\ \cup\ \dots\ \cup\ \hatcal N_r\ \cup\ \tildecal N, & \qquad &
			\hatcal F                                                                        & :=     & \hatcal F_0\ \cup\ \dots\ \cup\ \hatcal F_r\ \cup\
			\barcal F.
		\end{array}
	\]
	Let \( \hat N \) and \( \hat F \) be the free lattices over \(
	\hatcal N \) and \( \hatcal F \) respectively. In particular, we have
	\( \rank(\hat N) = 2n + (m-1)(r+1) \) and \( \rank(\hat F) = n + mr
	\). According to the type of \( P \), we define a map \( \gamma
	\colon \hat N \to \N \) as follows:

	\begin{center}

		\begin{longtable}{ll}

			(ee) & \(
			\begin{array}[t]{lcl}
				\hat e_{ij} & \mapsto &
				\begin{cases}
					e_i^+ - e_{i1},    & j = 1               \\
					e_{ij-1} - e_{ij}, & 2 \leq j \leq n_i-1
				\end{cases} \\[12pt]
				\hat u_{ij} & \mapsto &
				\begin{cases}
					u_i^+ - u_{i1},    & j = 1               \\
					u_{ij-1} - u_{ij}, & 2 \leq j \leq n_i-1
				\end{cases} \\[12pt]
				\tilde e_i  & \mapsto & e_i^+ - e_i^-    \\[3pt]
				\tilde u    & \mapsto & u^+ - u^-,
			\end{array} \) \\\\

			(pe) & \(
			\begin{array}[t]{lcl}
				\hat e_{ij} & \mapsto &
				\begin{cases}
					e_{ij-1} - e_{ij},  & 1 \leq j \leq n_i-1 \\
					e_{in_i-1} - e_i^-, & j = n_i
				\end{cases} \\[12pt]
				\hat u_{ij} & \mapsto &
				\begin{cases}
					u_{ij-1} - u_{ij}, & 1 \leq j \leq n_i-1 \\
					u_{in_i-1} - u^-,  & j = n_i
				\end{cases}
			\end{array} \)             \\\\

			(ep) & \(
			\begin{array}[t]{lcl}
				\hat e_{ij} & \mapsto &
				\begin{cases}
					e_i^+ - e_{i1},    & j = 1,            \\
					e_{ij-1} - e_{ij}, & 2 \leq j \leq n_i \\
				\end{cases} \\[12pt]
				\hat u_{ij} & \mapsto &
				\begin{cases}
					u^+ - u_{i1},      & j = 1,            \\
					u_{ij-1} - u_{ij}, & 2 \leq j \leq n_i \\
				\end{cases}
			\end{array} \)             \\\\

			(pp) & \(
			\begin{array}[t]{lcl}
				\hat e_{ij} & \mapsto & e_{ij-1} - e_{ij},                  \\[5pt]
				\hat u_{ij} & \mapsto & u_{ij-1} - u_{ij},                  \\[5pt]
				\tilde u_i  & \mapsto & u_{in_i} - u_{i-1,n_{i-1}}          \\[5pt]
				\tilde e    & \mapsto & (e_{00} + e_{10} + \dots + e_{r0}).
			\end{array} \)

		\end{longtable}

	\end{center}
	Then \( \gamma \) is a kernel of \( \alpha \colon \N \to N \). Next, we define
	a map \( \delta \colon \hat F \to \F \) as follows:
	\begin{center}
		\begin{longtable}{ll}
			(ee) & \(
			\begin{array}[t]{lcl}
				\hat f_{ij} & \mapsto & f_{ij}^+ - f_{ij}^-
			\end{array} \)                \\\\

			(pe) & \(
			\begin{array}[t]{lcl}
				\hat f_{ij} & \mapsto & f_{ij}^+ - f_{ij}^-     \\[5pt]
				\hat f_i^-  & \mapsto & f_{i-1,0}^- - f_{i,0}^-
			\end{array} \)     \\\\

			(ep) & \(
			\begin{array}[t]{lcl}
				\hat f_{ij} & \mapsto & f_{ij}^+ - f_{ij}^-
				\\[5pt]
				\hat f_i^+  & \mapsto & f_{i-1,n_{i-1}+1}^+ - f_{i,n_i+1}^+
			\end{array} \) \\\\

			(pp) & \(
			\begin{array}[t]{lcl}
				\hat f_{ij} & \mapsto & f_{ij}^+ - f_{ij}^-
				\\[5pt]
				\hat f_i^-  & \mapsto & f_{i-1,0}^- - f_{i,0}^-
				\\[5pt]
				\hat f_i^+  & \mapsto & f_{i-1,n_{i-1}+1}^+ - f_{i,n_i+1}^+.
			\end{array} \)             \\\\
		\end{longtable}
	\end{center}
	Then \( \delta \) is a kernel of \( \beta \colon \F \to F \). If \( P \) is of
	type (ee), set \( \tilde e_0 := -(\tilde e_1 + \dots + \tilde e_r) \). Then the
	induced map \( \hat P \colon \hat F \to \hat N \) is given as follows.
	\begin{center}

		\begin{longtable}{ll}
			(ee) & \(
			\begin{array}[t]{lcl}
				\hat f_{ij} & \mapsto &
				\begin{cases}
					l_{ij} \hat e_{ij} + d_{ij} \hat u_{ij}, & 1 \leq j \leq n_i-1 \\
					\displaystyle l_{in_i} \left( \tilde e_i - \sum_{k=0}^{n_i-1} \hat e_{ik}
					\right) + d_{in_i} \left( \tilde u - \sum_{k=1}^{n_i-1}
					\hat u_{ik} \right),                     & j = n_i,
				\end{cases}
			\end{array} \)                                                            \\\\

			(pe) & \(
			\begin{array}[t]{lc>{\displaystyle}l}
				\hat f_{ij} & \mapsto & l_{ij} \hat e_{ij} + d_{ij} \hat u_{ij}
				\\[5pt]
				\hat f_i^-  & \mapsto & \left( \sum_{k=1}^{n_{i-1}} \hat
				u_{i-1,k} \right) - \left( \sum_{k=1}^{n_i} \hat
				u_{ik} \right)
			\end{array} \) \\\\

			(ep) & \(
			\begin{array}[t]{lc>{\displaystyle}l}
				\hat f_{ij} & \mapsto & l_{ij} \hat e_{ij} + d_{ij} \hat u_{ij}
				\\[5pt]
				\hat f_i^+  & \mapsto & \left( \sum_{k=1}^{n_{i-1}} \hat
				u_{i-1,k} \right) - \left( \sum_{k=1}^{n_i} \hat u_{ik} \right)
			\end{array} \) \\\\

			(pp) & \(
			\begin{array}[t]{lc>{\displaystyle}l}
				\hat f_{ij} & \mapsto & l_{ij} \hat e_{ij} + d_{ij} \hat u_{ij}
				\\[5pt]
				\hat f_i^-  & \mapsto & \left( \sum_{k=1}^{n_{i-1}} \hat u_{i-1,k}
				\right) - \left( \sum_{k=1}^{n_i} \hat u_{ik} \right) - \tilde u_i \\
				\hat f_i^+  & \mapsto & \tilde u_i.
			\end{array} \)
		\end{longtable}
	\end{center}
\end{construction}

\begin{proof}
	In all cases, we can verify that \( \alpha \circ \gamma = 0 \) and \( \beta
	\circ \delta = 0 \). Moreover, \( \gamma \) and \( \delta \) are both
	injective and we have
	\begin{align*}
		\rank(\hat N) & = 2n - (m-1)(r+1) = \rank(\N) - \rank(N), \\
		\rank(\hat F) & = n + mr = \rank(\F) - \rank(F).
	\end{align*}
	This shows that \( \gamma \) and \( \delta \) are kernels of \( \alpha \)
	and \( \beta \) respectively. Furthermore, direct computations verify that \(
	\gamma \circ \hat P = \P \circ \delta \) holds in all cases. This shows that
	\( \hat P \) is the induced map between the kernels.
\end{proof}

\begin{example}
	\label{exm:running_example_phat}
	Continuing Example~\ref{exm:running_example_alpha_beta}, we get \(
	\hatcal N = \{\hat e_{01}, \hat u_{01}, \tilde e_1, \tilde e_2,
	\tilde u \} \) and \( \hatcal F = \{\hat f_{01}, \hat f_{02},
	\hat f_{11}, \hat f_{21} \} \). We obtain the matrix
	representations
	\[
		{\tiny
				\gamma=
				\begin{bmatrix}
					-1 & 0  & 1  & 0  & 0  \\
					-1 & 0  & 0  & 1  & 0  \\
					0  & 1  & 0  & 0  & 1  \\
					-1 & 0  & 0  & 0  & 0  \\
					0  & -1 & 0  & 0  & 0  \\
					0  & 0  & -1 & 0  & 0  \\
					0  & 0  & 0  & -1 & 0  \\
					0  & 0  & 0  & 0  & -1
				\end{bmatrix}, \quad
				\delta=
				\begin{bmatrix}
					1  & 0  & 0  & 0  \\
					0  & 0  & 1  & 0  \\
					0  & 0  & 0  & 1  \\
					-1 & 0  & 0  & 0  \\
					0  & 1  & 0  & 0  \\
					0  & -1 & 0  & 0  \\
					0  & 0  & -1 & 0  \\
					0  & 0  & 0  & -1
				\end{bmatrix},\quad
				\hat P=
				\begin{bmatrix}
					l_{01} & -l_{02} & 0      & 0      \\
					d_{01} & -d_{02} & 0      & 0      \\
					0      & -l_{02} & l_{11} & 0      \\
					0      & -l_{02} & 0      & l_{21} \\
					0      & d_{02}  & d_{11} & d_{21}
				\end{bmatrix}.
			}
	\]

\end{example}

\begin{remark}
	\label{rem:phat_matrix_representation}
	We give the general matrix representation of \( \hat P \) from
	Construction~\ref{cns:phat_cstar_surface} for the case (ee).
	Let \( \hat F_i \) and \( \hat N_i \) be the free lattices over \( \hatcal F_i \)
	and \( \hatcal N_i \) respectively. Let \( \tilde N \) be the free lattice
	over \( \tildecal N \). We define the lattice maps
	\[
		\begin{array}{lclcl}
			\hat P_i \colon \hat F_i \to \hat N_i,   & \qquad  &
			\hat f_{ij}                              & \mapsto &
			\begin{cases}
				l_{ij} \hat e_{ij} + d_{ij} \hat u_{ij},       & 1 \leq j \leq n_i-1 \\
				\displaystyle -l_{in_i} \left( \sum_{k=1}^{n_{i}-1} \hat e_{ik} \right) - d_{in_i}
				\left( \sum_{k=1}^{n_i-1} \hat u_{ik} \right), & j = n_i,
			\end{cases} \\[20pt]
			\tilde P_i \colon \hat F_i \to \tilde N, & \qquad  &
			\hat f_{ij}                              & \mapsto &
			\begin{cases}
				0,                                       & 1 \leq j \leq n_i-1 \\
				l_{in_i} \tilde e_i + d_{in_i} \tilde u, & j = n_i.
			\end{cases}.
		\end{array}
	\]
	Then we have \( \hat P(\hat f_{ij}) = \hat P_i(\hat f_{ij}) + \tilde
	P_i(\hat f_{ij}) \). We obtain the matrix representations:
	\[
		\begin{array}{lclclcl}
			\tilde P_0                                                              & =      &
			\begin{bmatrix}
				0      & \dots & 0      & -l_{0n_0} \\
				0      & \dots & 0      & -l_{0n_0} \\
				\vdots &       & \vdots & \vdots    \\
				0      & \dots & 0      & -l_{0n_0} \\
				0      & \dots & 0      & d_{0n_0}
			\end{bmatrix},               & \qquad &
			\tilde P_i                                                              & =      &
			\begin{bmatrix}
				0      & \dots & 0      & 0        \\
				\vdots &       & \vdots & \vdots   \\
				0      & \dots & 0      & l_{in_i} \\
				\vdots &       & \vdots & \vdots   \\
				0      & \dots & 0      & 0        \\
				0      & \dots & 0      & d_{in_i}
			\end{bmatrix} \quad (i \geq 1),                            \\\\

			\hat P_i                                                                & =      &
			\begin{bmatrix}
				l_{i1} & 0      & \dots  & 0          & -l_{in_i} \\
				d_{i1} & 0      & \dots  & 0          & -d_{in_i} \\
				0      & l_{i2} &        & 0          & -l_{in_i} \\
				0      & d_{i2} &        & 0          & -d_{in_i} \\
				\vdots &        & \ddots &            & \vdots    \\
				0      & 0      & \dots  & l_{in_i-1} & -l_{in_i} \\
				0      & 0      & \dots  & d_{in_i-1} & -d_{in_i}
			\end{bmatrix}, & \qquad &
			\hat P                                                                  & =      &
			\begin{bmatrix}
				\hat P_0   & 0          & \dots  & 0          \\
				0          & \hat P_1   & \dots  & 0          \\
				\vdots     &            & \ddots &            \\
				0          & 0          &        & \hat P_r   \\
				\tilde P_0 & \tilde P_1 & \dots  & \tilde P_r
			\end{bmatrix}.
		\end{array}
	\]
\end{remark}

\section{Maximal minors of \texorpdfstring{\( \hat P \)}{}}
\label{sec:maximal_minors_phat}

In this section, we consider the set $M(\hat P)$ of maximal minors
arising from the matrix representation of $\hat P$ given in the
preceding section. A series of reduction steps turns $M(\hat P)$ into
a smaller set having the same greatest common divisor, which finally
enables us to compare with $M'(P)$. Whereas the cases (pe), (ep) and
(pp) are are the simpler ones; see
Proposition~\ref{prp:maximal_nonzero_minors_phat_pe_pe_pp}, the case
(ee) takes up the second half of the section, see
Proposition~\ref{prp:gcd_minors_prime_union_minors_prime_prime_succ}.

\begin{definition}
	\label{def:maximal_minors_Phat_cstar_surface}
	Let a subset \( A \subseteq \hatcal N \) with \( |A| = |\hatcal F| \) be given. Then we
	have a sublattice \( \hat N_A := \bigoplus_{x \in A} \ZZ \cdot x \subseteq \hat N \) and
	an induced map \( \hat P_A \colon \hat F \to \hat N_A \) as in the commutative diagram
	\[
		\begin{tikzcd}
			\hat F \ar[r, "\hat P"] \ar[dr, "\hat P_A" {below,yshift=-0.5ex,near
			start}] & \hat N                           \\
			        & \hat N_A \ar[u, hookrightarrow].
		\end{tikzcd}
	\]
	We call \( |\det(\hat P_A)| \in \ZZ \) the \emph{maximal minor of \(
		\hat P \) associated to \( A \)}. The set of all maximal minors of \(
	\hat P \) is defined as
	\[
		M(\hat P) := \{ |\det(\hat P_A)| \ ;\ A \subseteq \hatcal N,\ |A| = |\hatcal F|
		\}.
	\]
\end{definition}

\begin{construction}
	\label{cns:reduced_maximal_minor_phat}
	Let \( A \subseteq \hatcal N \) with \( |A| = |\hatcal F| \). We define
	\begin{align*}
		\hatcal N_A^{\sing}\ :=\  & \{\hat e_{ij}\ ;\ \hat e_{ij} \in A \text{ and
		} \hat u_{ij} \notin A \}\ \cup\ \{\hat u_{ij}\ ;\ \hat e_{ij} \notin
		A \text{ and } \hat u_{ij} \in A \}\ \subseteq\ A,                         \\
		\hatcal F_A^{\sing}\ :=\  & \{\hat f_{ij}\ ;\ \hat e_{ij} \in A \text{ and
		} \hat u_{ij} \notin A \}\ \cup\ \{\hat f_{ij}\ ;\ \hat e_{ij} \notin
		A \text{ and } \hat u_{ij} \in A \}\ \subseteq\ \hatcal F,
	\end{align*}
	\[
		\hatcal N_A^{\red}\ :=\ A\ \backslash\ \hatcal N_A^{\sing}, \qquad \qquad
		\hatcal F_A^{\red}\ :=\ \hatcal F\ \backslash\ \hatcal F_A^{\sing}.
	\]
	Note that we have \( |\hatcal N_A^{\red}| = |\hatcal F_A^{\red}| \).
	Let \( \hat N_A^{\red} \) and \( \hat F_A^{\red} \) be the free
	lattices over \( \hatcal N_A^{\red} \) and \( \hatcal F_A^{\red} \)
	respectively. We obtain an induced map \( \hat P_A^{\red} \) as in
	the commutative diagram
	\[
		\begin{tikzcd}
			\hat F \ar[r, "\hat P_A"] & \hat N_A \\
			\hat F_A^{\red} \ar[u, hookrightarrow] \ar[r, "\hat P_A^{\red}"
			{below}]                  &
			\hat N_A^{\red} \ar[u, hookrightarrow].
		\end{tikzcd}
	\]
	We call \( |\det(\hat P_A^{\red})| \) the \emph{reduced minor of \(
		\hat P \) associated to \( A \)}. We define the set of reduced minors
	of \( \hat P \) as
	\[
		M^{\red}(\hat P) := \{ |\det(\hat P_A^{\red})| \ ;\ A \subseteq \hatcal N,\
		|A| = |\hatcal F|\}.
	\]
\end{construction}

\begin{example}
	\label{exm:running_example_maximal_minors_phat}
	The matrix \( \hat P \)
	from Example~\ref{exm:running_example_phat} has five maximal
	minors and four reduced minors:
	\begin{align*}
		M(\hat P)        & = \{d_{01}|\hat\mu|,\quad l_{01}|\hat\mu|,\quad
		d_{11}l_{21}|\hat\nu(0,1)|,\quad
		l_{11}d_{21}|\hat\nu(0,1)|,\quad
		l_{11}l_{21}|\hat\nu(0,1)| \}                                      \\
		M^{\red}(\hat P) & = \{|\hat\mu|,\quad
		d_{11}l_{21}|\hat\nu(0,1)|,\quad
		l_{11}d_{21}|\hat\nu(0,1)|,\quad
		l_{11}l_{21}|\hat\nu(0,1)| \}.
	\end{align*}
\end{example}

\begin{proposition}
	\label{prp:reduced_maximal_minor_phat}

	Let \( A \subseteq \hatcal N \) with \( |A| = |\hatcal F| \). Then we
	have
	\[
		\det(\hat P_A) = \det(\hat P_A^{\red}) \prod_{\hat f_{ij} \in \hatcal
		F_A^{\sing}} x_{ij}, \quad
		\text{where} \quad x_{ij} :=
		\begin{cases}
			l_{ij}, & \hat e_{ij} \in A, \\
			d_{ij}, & \hat u_{ij} \in A.
		\end{cases}
	\]
	Moreover, \( \gcd(M(\hat P)) = \gcd(M^{\red}(\hat P)) \) holds.
\end{proposition}

\begin{proof}
	If \( \hat f_{ij} \in \hatcal F_A^{\sing} \), we have
	\[
		\hat P_A(\hat f_{ij}) =
		\begin{cases}
			l_{ij}\hat e_{ij}, & \hat e_{ij} \in A, \\
			d_{ij}\hat u_{ij}, & \hat u_{ij} \in A.
		\end{cases}
	\]
	In other words, the matrix representation of \( \hat P_A \) has a
	column with the single entry \( x_{ij} \) and zeroes elsewhere. Doing
	cofactor expansion by all these columns amounts to passing from \(
	\det(\hat P_A) \) to \( \det(\hat P_A^{\red}) \). This shows the
	first claim. The second one then follows from \( \gcd(l_{ij}, d_{ij})
	= 1 \).
\end{proof}

\begin{definition}
	\label{def:maximal_minors_phat_leaf_configuration}
	Set \(\cal L := \{(i,j)\ ;\ i = 0, \dots, r,\ j = 1, \dots, n_i' \} \).
	For \( A \subseteq \hat{\cal N} \) with \( |A| = |\hatcal F| \), we define
	\[
		L(A)\ :=\ \{(i,j)\ ;\ \hat e_{ij} \in A \text{ and } \hat u_{ij} \in A
		\}\ \subseteq\ \cal L.
	\]
\end{definition}

\begin{lemma}
	\label{lem:maximal_minors_phat_vanishing}
	Let \( A \subseteq \hat{\cal N} \) with \( |A| = |\hatcal F| \).
	\begin{enumerate}[\normalfont(i)]
		\item If \( \hat e_{ij} \notin A \) and \( \hat u_{ij} \notin A \) for some
		      \( i = 0, \dots, r \) and \( j = 1, \dots, n_i' \), we have \(
		      \det(\hat P_A) = 0 \).
		\item If \( (i,j_0), (i,j_1) \in L(A) \) for some \( i = 0, \dots, r \) and
		      \( 1 \leq j_0 < j_1 \leq n_i' \), we have \( \det(\hat P_A) = 0 \).
	\end{enumerate}
\end{lemma}

\begin{proof}
	For (i), we have \( \hat P_A(\hat f_{ij}) = 0 \), hence \( \det(\hat P_A) = 0
	\). We show (ii). Consider first the case (ee). Then we have \( n_i' = n_i
	- 1 \). Set \( \hatcal N_{A,i} := A \cap \hatcal N_i \) and \(
	\hat N_{A,i} := \hat N_A \cap \hat N_i \). By (i),
	we may assume that \( \hat e_{ij} \in \hatcal N_{A,i} \) or \( \hat u_{ij}
	\in \hatcal N_{A,i} \) holds for all \( 1 \leq j \leq n_i-1 \). Since \(
	(i,j_0), (i,j_1) \in L(A) \), we thus have \( |\hatcal N_{A,i}| > n_i \).
	Consider the map
	\( \hat P_{A,i} \colon \hat N_{A,i} \to \hat F_i \). In the
	matrix representation from Remark \ref{rem:phat_matrix_representation}, we have
	\[
		\hat P_A =
		\begin{bmatrix}
			\ast & 0                                                        & 0
			\\
			\ast & \boxed{\begin{array}{c} \hat P_{A,i} \quad 0\end{array}} &
			0
			\\
			\ast & \ast                                                     & \ast
		\end{bmatrix},
	\]
	where the outlined box is a square \( |\hatcal N_{A,i}| \times
	|\hatcal N_{A,i}| \)-matrix. Since the determinant of the outlined
	box vanishes, also \( \det(\hat P_A) = 0 \).

	Now let \( P \) be of type (pe), (ep) or (pp). Then \( n_i' = n_i \).
	We define the set
	\[
		\barcal N_A\ :=\ \{ \hat u_{ij}\ ;\ (i,j) \in L(A) \}\ \subseteq\
		\hatcal N_A^{\red}.
	\]
	Writing \( \bar N_A \) for the free lattice over \( \barcal N_A \)
	and \( \bar F \) for the free lattice over \( \barcal F \), we obtain
	an induced map \( \bar P_A \colon \bar F \to \bar N_A \) as in the
	commutative diagram
	\[
		\begin{tikzcd}
			\hat F_A^{\red} \ar[r, "\hat P_A^{\red}"]                & \hat N_A^{\red} \\
			\bar F \ar[u, hookrightarrow] \ar[r, "\bar P_A" {below}] & \bar N_A
			\ar[u, hookrightarrow].
		\end{tikzcd}
	\]
	Note that if \( (i,j) \in L(A) \), we have \( (\hat
	P_A^{\red})^*(\hat e_{ij}^*) = l_{ij} \hat f_{ij}^* \). That means,
	\( \hat P_A^{\red} \) contains a row with a single entry \( l_{ij} \)
	and zeroes elsewhere. Doing cofactor expansion, we arrive at
	\[
		\det(\hat P_A^{\red}) = \det(\bar P_A) \prod_{(i,j) \in L(A)} l_{ij}.
	\]
	But since \( (i,j_0), (i,j_1) \in L(A) \), we have \( \bar P_A^*(\hat
	u_{ij_0}^*) = \bar P_A^*(\hat u_{ij_1}^*) \), i.e. \( \bar P_A \)
	contains two equal rows. Hence we have \( \det(\hat P_A) = \det(\hat
	P_A^{\red}) = \det(\bar P_A) = 0 \).
\end{proof}

\begin{proposition}
	\label{prp:maximal_nonzero_minors_phat_pe_pe_pp}
	Let \( P \) be of type (pe), (ep) or (pp). Let \( A \subseteq \hatcal N
	\) with \( |A| = |\hatcal F| \) such that \( \det(\hat P_A) \neq 0 \). Then
	we have \( |L(A)| = r \). Furthermore, there exists an \( i_1 = 0, \dots, r
	\) and \( j_i = 1, \dots, n_i \) for all \( i \neq i_1 \) such that
	\[
		|\det(\hat P_A^{\red})| = \prod_{i \neq i_1} l_{ij_i}.
	\]
	In particular, we have \( \gcd(M^{\red}(\hat P)) = \gcd(M'(P)) \).
\end{proposition}

\begin{proof}

	By Lemma \ref{lem:maximal_minors_phat_vanishing}~(i), we have \( \hat
	e_{ij} \in A \) or \( \hat u_{ij} \in A \) for all \( i \) and \( j
	\). By Lemma~\ref{lem:maximal_minors_phat_vanishing}~(ii), for each
	\( i \) there is at most one \( j \) with \( \hat e_{ij} \in A \) and
	\( \hat u_{ij} \in A \). Writing \( \pi_1 \colon \ZZ \times \ZZ \to
	\ZZ \) for the projection onto the first coordinate, this implies
	that \( |\pi_1(L(A))| = |L(A)| \). Together, we obtain
	\[
		|A \cap \hatcal N_i| =
		\begin{cases}
			n_i,     & i \notin \pi_1(L(A)), \\
			n_i + 1, & i \in \pi_1(L(A)).
		\end{cases}
	\]
	Thus, we have \( |A| = n + |L(A)| + |A \cap \tildecal N| \). Recall
	that for the cases (pe) and (ep), we have \( \tildecal N = \emptyset
	\) and \( |A| = |\hatcal F| = n + r \), hence \( |L(A)| = r \). For
	(pp), we must have \( \tilde u_i \in A \) for all \( i = 1, \dots, r
	\), since otherwise \( \hat P_A(\hat f_i^+) = 0 \). Hence \( |A \cap
	\tildecal N| = r \). Since in this case, \( |A| = n+2r \), we also
	arrive at \( |L(A)| = r \). This implies that we have \( L(A) =
	\{(i,j_i)\ ;\ i \neq i_1\} \) for some \( i_1 = 0, \dots, r \) and \(
	j_i = 1, \dots, n_i \). Following the proof of Lemma
	\ref{lem:maximal_minors_phat_vanishing} (ii), we proceed to do
	cofactor expansion and see that \( |\det(\bar P_A)| = 1 \). This
	proves the claim. The supplement follows directly from the Definition
	of \( M'(P) \).
\end{proof}

The preceding proposition settles the discussion of maximal minors of
\( \hat P \) for the cases (pe), (pe) and (pp). The remainder of this
section is devoted to the case (ee).

\begin{lemma}
	\label{rem:maximal_minors_phat_count_tilde}
	Let \( P \) be of type (ee). Let \( A \subseteq \hat{\cal N} \) with \( |A| =
	|\hatcal F| = n \) and \( \det(\hat P_A) \neq 0 \). Then we have
	\[
		|L(A)| = (r+1) - |A \cap \tilde{\cal N}|.
	\]
	In particular, \( 0 \leq |L(A)| \leq r+1 \).
\end{lemma}

\begin{proof}
	Let \( \pi_1 \colon \ZZ \times \ZZ \to \ZZ \) be the projection onto the first
	coordinate. Lemma \ref{lem:maximal_minors_phat_vanishing}
	(ii) implies \( |\pi_1(L(A))| = |L(A)|. \) Furthermore, we have
	\[
		|A \cap \hat{\cal N}_i| =
		\begin{cases}
			n_i-1, & i \notin \pi_1(L(A)) \\
			n_i,   & i \in \pi_1(L(A))
		\end{cases}.
	\]
	Since, \( |A| = n \), this implies the claim.
\end{proof}

\begin{definition}
	\label{def:maximal_minors_reduced_k}
	Let \( P \) be of type (ee). For \( k = 0, \dots, r+1 \), we define
	\[
		M^{\red}_k(\hat P) := \{ |\det(\hat P^{\red}_A)|\ ;\ |L(A)| = k \}.
	\]
\end{definition}

\begin{remark}
	\label{rem:reduced_maximal_minors_zero_eq_mu_hat}
	Let \( A \subseteq \hatcal F \) with \( |A| = n \) and \( L(A) = \emptyset
	\). Lemma \ref{rem:maximal_minors_phat_count_tilde} implies that \( A \cap
	\tildecal N = \tildecal N \). We obtain \( \det(\hat P_A^{\red}) = \hat \mu
	\). Hence \( M^{\red}_0(\hat P) = \{|\hat \mu|\} \).
\end{remark}

\begin{proposition}
	\label{prp:maximal_minors_phat_reduced_phat_types_cstar_surface_ee}
	Let \( P \) be of type (ee). Let \( A \subseteq \hat{\cal N} \) with \( |A|=n
	\) such that \( |\hat P_A^{\red}| \neq 0 \). Write \( \pi_1 \colon \ZZ \times \ZZ
	\to \ZZ \) for the projection onto the first coordinate.
	\begin{enumerate}[\normalfont(i)]
		\item Assume \( \tilde u \notin A \). Then for all \( i = 1, \dots, r \),
		      we have \( \tilde e_i \in A \) or \( i \in \pi_1(L(A)) \). In this
		      case,
		      \[
			      |\det(\hat P_A^{\red})| = \left( \prod_{(i,j) \in L(A)}
			      |\hat \nu(i,j)| \right) \left( \prod_{\substack{0 \leq i \leq r \\ i \notin
					      \pi_1(L(A))}} l_{in_i}\right).
		      \]
		\item Assume \( \tilde u \in A \). Then there exists at most one \( i_1 =
		      1, \dots, r \) with \( \tilde e_{i_1} \notin A \) and \( i_1 \notin
		      \pi_1(L(A)) \). If there exists such an \( i_1 \), we have
		      \[
			      |\det(\hat P_A^{\red})| = \left| d_{i_1n_{i_1}} \left( \prod_{(i,j) \in
					      L(A)} \hat \nu(i,j) \right) \left( \prod_{\substack{0 \leq i
					      \leq r \\ i \notin \pi_1(L(A)) \cup \{i_1\}}} l_{in_i}
			      \right)\right|.
		      \]
		      If there exists no such \( i_1 \), we have
		      \[
			      |\det(\hat P_A^{\red})| = \left| \sum_{\substack{0 \leq i' \leq r \\ i'
					      \notin \pi_1(L(A))}} \pm d_{i'n_{i'}} \left( \prod_{(i,j) \in
					      L(A)} \hat \nu(i,j) \right) \left( \prod_{\substack{0 \leq i
					      \leq r \\ i \notin \pi_1(L(A)) \cup \{i'\}}}
			      l_{in_i}
			      \right) \right|.
		      \]
	\end{enumerate}

	\begin{proof}
		For part (i), assume that there is some \( i = 1, \dots, r \) such that \(
		\tilde e_i \notin A \) and \( i \notin \pi_1(L(A)) \). This implies \( \hat
		e_{ij} \notin \hatcal N_A^{\red} \) and \( \hat u_{ij} \notin \hatcal
		N_A^{\red} \) for all \( j = 1, \dots, n_i-1 \). Since also \( \tilde u \notin
		A \), we obtain \( \hat P_A^{\red}(\hat f_{in_i}) = 0 \). Hence \( \det(\hat
		P_A^{\red}) = 0 \), a contradiction. It follows that if \( i \notin \pi_1(L(A))
		\), either \( i = 0 \) or \( \tilde e_i \in A \). The formula for \( \det(\hat
		P_A^{\red}) \) now follows from cofactor expansion.

		For part (ii), assume that there exist \( 1 \leq i_0 < i_1 \leq r \)
		such that \( \tilde e_{i_0}, \tilde e_{i_1} \notin A \) and \( i_0,
		i_1 \notin \pi_1(L(A)) \). We obtain \( \hat P_A^{\red}(\hat
		f_{i_0n_{i_0}}) = d_{i_0n_{i_0}} \tilde u \) and \( \hat
		P_A^{\red}(\hat f_{i_1n_{i_1}}) = d_{i_1n_{i_1}} \tilde u \). Hence
		\( \det(\hat P_A^{\red}) = 0 \), a contradiction. The formulas for \(
		\det(\hat P_A^{\red}) \) again follow from cofactor expansion.
	\end{proof}
\end{proposition}

\begin{definition}
	\label{def:maximal_minors_prime_phat}
	Let \( P \) be of type (ee). Let \( \pi_1 \colon \ZZ \times \ZZ \to \ZZ \) be
	the projection onto the first coordinate. We call a subset \( L \subseteq \cal
	L \) \emph{valid}, if \( |\pi_1(L)| = |L| \). For \( k = 1, \dots, r \), we
	define
	\[
		M'_k(\hat P) := \left \{ \left( \prod_{(i,j) \in L} |\hat \nu(i,j)| \right)
		\left( \prod_{\substack{0 \leq i \leq r \\ i \notin \pi_1(L) \cup
				\{i_1\}}} l_{in_i} \right)\ ;\
		\begin{aligned}
			 & L \subseteq \cal L \text{ valid},\ |L| = k, \\
			 & 0 \leq i_1 \leq r,\ i_1 \notin \pi_1(L),
		\end{aligned}
		\right\}.
	\]
\end{definition}

\begin{example}
	\label{exm:running_example_maximal_minors_phat_prime}
	Continuing Example~\ref{exm:running_example_maximal_minors_phat},
	we have
	\[
		M_0^{\red}(\hat P)\ =\ \{|\hat\mu|\}, \quad M_1^{\red}(\hat
		P)\ =\
		\{d_{11}l_{21}|\hat\nu(0,1)|,\ l_{11}d_{21}|\hat\nu(0,1)|,\
		l_{11}l_{21}|\hat\nu(0,1)|\}
	\]
	\[
		M_1'(\hat P)\ =\ \{l_{11}|\hat\nu(0,1)|,\
		l_{21}|\hat\nu(0,1)|\},\qquad M_2'(\hat P) = M_2^{\red}(\hat
		P) = M_3^{\red}(\hat P) = \emptyset.
	\]
\end{example}

\begin{proposition}
	\label{prp:gcd_maximal_minors_phat_eq_gcd_maximal_minors_prime}
	Let \( P \) be of type (ee). We have
	\begin{enumerate}[\normalfont(i)]
		\item \(\gcd(M^{\red}_k(\hat P)) = \gcd(M'_k(\hat P)) \) for all \( k =
		      1, \dots, r \),
		\item \( \gcd(M^{\red}_{r+1}(\hat P) \cup M'_r(\hat P)) = \gcd(M'_r(\hat P)) \),
		\item \( \gcd\left(\bigcup_{k=1}^{r+1} M_k^{\red}(\hat P)\right) =
		      \gcd\left(\bigcup_{k=1}^{r} M_k'(\hat P)\right) \).
	\end{enumerate}
\end{proposition}

\begin{proof}
	We show (i). Proposition
	\ref{prp:maximal_minors_phat_reduced_phat_types_cstar_surface_ee} implies
	that every element of \( M_k^{\red}(\hat P) \) is a \( \ZZ \)-linear
	combination of elements of \( M'_k(\hat P) \). This shows that \( \gcd(M'_k(\hat
	P)) \) divides \( \gcd(M_k^{\red}(\hat P)) \). For the converse, it suffices to show
	that \( \gcd(M_k^{\red}(\hat P)) \mid x \) holds for all \( x \in M'_k(\hat
	P) \). So, let
	\[
		x = \left( \prod_{(i,j) \in L} \nu(i,j) \right)
		\left( \prod_{\substack{0 \leq i \leq r \\ i \notin \pi_1(L) \cup
				\{i_1\}}} l_{in_i} \right) \in M'_k(\hat P)
	\]
	be arbitraty, where \( L \subseteq \cal L \) is a valid subset with
	\( |L| = k \) and \( 0 \leq i_1 \leq r \) with \( i_1 \notin \pi_1(L)
	\).

	\emph{Case 1: \( i_1 \neq 0 \) and \( 0 \in \pi_1(L) \) }. Choose a subset
	\( A \subseteq \hat{\cal N} \) with \( |A| = n \) such that \( L(A) = L \)
	and
	\[
		A \cap \tilde{\cal N} = \{\tilde e_i\ ;\ 1 \leq i \leq r,\ i \notin
		\pi_1(L) \}.
	\]
	Since \( 0 \in \pi_1(L) \), we have \( |A \cap \tilde{\cal N}| = r -
	(|L| - 1) = (r+1) - |L| \). This means we can choose \( A \) such
	that \( \det(\hat P^{\red}_A) \neq 0 \). Now set
	\[
		A' := (A \backslash \{\tilde e_{i_1}\}) \cup \{\tilde u\}.
	\]
	Proposition
	\ref{prp:maximal_minors_phat_reduced_phat_types_cstar_surface_ee}
	implies that \( \det(\hat P_A^{\red}) = l_{i_1n_{i_1}} x \) and \(
	\det(\hat P_{A'}^{\red}) = d_{i_1n_{i_1}} x \). Hence we have
	\[
		\gcd(M_k^{\red}(\hat P)) \mid \gcd(l_{i_1n_{i_1}} x, d_{i_1n_{i_1}} x) =
		x.
	\]

	\emph{Case 2: \( i_1 \neq 0 \) and \( 0 \notin \pi_1(L) \) }. Since \( k
	\geq 1 \), we find some \( i_0 \in \pi_1(L) \). Now choose \( A \subseteq
	\hat{\cal N} \) with \( |A| = n \) such that \( L(A) = L \) and
	\[
		A \cap \tilde{\cal N} = \{\tilde e_i\ ;\ 1 \leq i \leq r,\ i \notin
		\pi_1(L) \} \cup \{ \tilde e_{i_0} \}.
	\]
	Again, we have \( |A \cap \tilde{\cal N}| = (r+1) - |L| \), hence we
	can pick \( A \) such that \( \det(\hat P_A^{\red}) \neq 0 \).
	Proceeding in the same way as in Case 1, we arrive at \(
	\gcd(M_k^{\red}(\hat P)) \mid x \).

	\emph{Case 3: \( i_1 = 0 \)}. As in Case 2, we can pick some \( i_0 \in
	\pi_1(L) \) as well as a subset \( A \subseteq \hat{\cal N} \) with
	\( |A| = n \) such that \( L(A) = L \) and
	\[
		A \cap \tilde{\cal N} = \{\tilde e_i\ ;\ 1 \leq i \leq r,\ i \notin
		\pi_1(L) \} \cup \{ \tilde e_{i_0} \}.
	\]
	For all \( 1 \leq i \leq r \) with \( i \notin \pi_1(L) \), set \(A_i
	:= (A \backslash \{\tilde e_i\}) \cup \{ \tilde u \}.\) Then
	Proposition
	\ref{prp:maximal_minors_phat_reduced_phat_types_cstar_surface_ee}
	implies that \( \det(P_A^{\red}) = l_{0n_{0}} x \) and
	\[
		\det(P_{A_{i_0}}^{\red}) = d_{0n_0} x + \sum_{\substack{1 \leq i' \leq r \\
				i' \notin \pi_1(L)}} \pm \det(P_{A_{i'}}^{\red}).
	\]
	This implies \(\gcd(M_k^{\red}(\hat P)) \mid \gcd(l_{0n_0} x,
	d_{0n_0} x) = x \).

	Part (ii) follows from the fact that every element of \(
	M^{\red}_{r+1}(\hat P) \) is an integer multiple of an element of \(
	M'_r(\hat P) \). Part (iii) is a consequence of (i) and (ii).
\end{proof}

\begin{definition}
	\label{def:maximal_minors_prime_prime_cstar_surface_ee}
	Let \( P \) be of type (ee). For \( k = 1, \dots, r \), we define the set
	\[
		M''_k(\hat P) := \left \{ \left( \prod_{(i,j) \in L} |\hat \nu(i,j)| \right)
		\left( \prod_{\substack{0 \leq i \leq r \\ i \notin \pi_1(L) \cup
				\{i_1\}}} l_{ij_i} \right)\ ;\
		\begin{aligned}
			 & L \subseteq \cal L \text{ valid},\ |L| = k \\
			 & 0 \leq i_1 \leq r,\ i_1 \notin \pi_1(L)    \\
			 & 1 \leq j_i \leq n_i \text{ for all } i
		\end{aligned}
		\right\}.
	\]
\end{definition}

\begin{lemma}
	\label{lem:gcd_lini_nu_hat_div_lij}
	Let \( i = 0, \dots, r \) and \( j_i = 1, \dots, n_i-1 \). Then we have
	\[
		\gcd(l_{in_i}, \hat \nu(i,j)) \mid l_{ij}.
	\]
\end{lemma}

\begin{proof}
	By definition we have \( \hat \nu(i,j) = l_{in_i} d_{ij} - l_{ij} d_{in_i} \).
	Since \( l_{in_i} \) and \( d_{in_i} \) are coprime, we find \( x,y \in \ZZ \)
	such that \( x l_{in_i} + y d_{in_i} = 1 \). Then we have
	\[
		(x l_{ij} + y d_{in_i}) l_{in_i} - y \hat \nu(i,j) = l_{ij}(x l_{in_i} +
		y d_{in_i}) = l_{ij}.
	\]
	This implies the claim.
\end{proof}

\begin{proposition}
	\label{prp:gcd_minors_prime_union_minors_prime_prime_succ}
	Let \( P \) be of type (ee). We have
	\begin{enumerate}[\normalfont(i)]
		\item \(\gcd(M'_k(\hat P) \cup M''_{k+1}(\hat P)) = \gcd(M''_k(\hat
		      P)) \) for all \( k = 1, \dots, r-1 \),
		\item \( \gcd\left(\bigcup_{k=1}^r M'_k(\hat P)\right) =
		      \gcd(M''_1(\hat P)) \).
	\end{enumerate}
\end{proposition}

\begin{proof}
	We show (i). Since \( M'_k(\hat P) \subseteq M''_k(\hat P) \) and elements of
	\( M''_{k+1}(\hat P) \) are \( \ZZ \)-linear combinations of elements of \(
	M''_k(\hat P) \), we have \( \gcd(M''_k(\hat P)) \mid \gcd(M'_k(\hat P) \cup
	M''_{k+1}(\hat P)) \). For the converse, let
	\[
		x = \left( \prod_{(i,j) \in L} \hat \nu(i,j) \right)
		\left( \prod_{\substack{0 \leq i \leq r \\ i \notin \pi_1(L) \cup
				\{i'\}}} l_{ij_i} \right) \in M''_k(\hat P),
	\]
	where \( L \subseteq \cal L \) is a valid subset with \( |L| = k \)
	and \( 0 \leq i' \leq r \) with \( i' \notin \pi_1(L) \) and \( 1
	\leq j_i \leq n_i \) for all \( i \notin \pi_1(L) \cup \{i'\} \). Let
	us write
	\[
		\{i_1, \dots, i_{r-k} \} := \{0, \dots, r\} \backslash (\pi_1(L) \cup
		\{i'\}).
	\]
	We define numbers
	\[
		\begin{array}{lcl}
			x_0     & :=     & l_{i_1n_{i_1}} l_{i_2n_{i_2}} \dots l_{i_{r-k}n_{i_{r-k}}} \\[4pt]
			x_1     & :=     & l_{i_1j_{i_1}} l_{i_2n_{i_2}} \dots l_{i_{r-k}n_{i_{r-k}}} \\[4pt]
			        & \vdots &                                                            \\[4pt]
			x_{r-k} & :=     & l_{i_1j_{i_1}} l_{i_2j_{i_2}} \dots l_{i_{r-k}j_{i_{r-k}}} \\
		\end{array}
	\]
	as well as
	\[
		\begin{array}{lcl}
			y_1     & :=     & \hat \nu(i_1, j_{i_1}) l_{i_2n_{i_2}}
			l_{i_3n_{i_3}} \dots l_{i_{r-k}n_{i_{r-k}}}                             \\[4pt]
			y_2     & :=     & l_{i_1j_{i_1}} \hat \nu(i_2, j_{i_2})
			l_{i_3n_{i_3}} \dots l_{i_{r-k}n_{i_{r-k}}}                             \\[4pt]
			        & \vdots                                                        \\[4pt]
			y_{r-k} & :=     & l_{i_1j_{i_1}} \dots l_{i_{r-k-1}j_{i_{r-k-1}}} \hat
			\nu(i_{r-k}, j_{i_{r-k}}).
		\end{array}
	\]
	Then Lemma \ref{lem:gcd_lini_nu_hat_div_lij} implies \(
	\gcd(x_{m-1},y_m) \mid x_m \) for all \( m = 1, \dots, r-k \). In
	particular, we obtain \( \gcd(x_0, y_1, \dots, y_{r-k}) \mid x_{r-k}
	\). Now set \( c := \prod_{(i,j) \in L} \hat \nu(i,j) \). Then we
	have \( cx_{r-k} = x \) as well as \( x_0c \in M'_k(\hat P) \) and \(
	y_mc \in M''_{k+1}(\hat P) \) for all \( m = 1, \dots r-k \).
	Together, we have
	\[
		\gcd(M'_k(\hat P) \cup M''_{k+1}(\hat P)) \mid \gcd(x_0c, y_1c, \dots,
		y_{r-k}c) \mid x_{r-k}c = x.
	\]
	Part (ii) follows from repeated application of (i), together with the
	fact that \( M'_r(\hat P) = M''_r(\hat P) \).
\end{proof}

\section{Proof of Theorem \ref{thm:picard_index_formula} and examples} \label{sec:proof_of_picard_index_formula}

In this section, we prove the formula for the Picard index of a \(
\KK^* \)-surface given in Theorem \ref{thm:picard_index_formula}. We
then give two examples where the formula fails: The first one is a
toric threefold, the second one is the \( D_8 \)-singular log del
Pezzo surface of Picard number one.

\begin{proposition}
	\label{prp:gcd_maximal_minors_P_eq_gcd_maximal_minors_Phat}

	Let \( P \) be a defining matrix and \( \hat P \) be as in
	Construction~\ref{cns:phat_cstar_surface}. Write \( M(P) \) and \(
	M(\hat P) \) for the set of maximal minors of \( P \) and \( \hat P
	\) respectively. Then we have
	\[
		\gcd(M(\hat P)) = \gcd(M(P)).
	\]
\end{proposition}

\begin{proof}
	By Construction~\ref{cns:reduced_maximal_minor_phat}, we have \( \gcd(M(\hat
	P)) = \gcd(M^{\red}(\hat P)) \). On the other hand, Proposition~\ref{prp:gcd_maximal_minors_P_eq_gcd_maximal_minors_prime}
	says that \( \gcd(M(P)) = \gcd(M'(P)) \). In the cases (pe), (ep) and (pp),
	Proposition~\ref{prp:maximal_nonzero_minors_phat_pe_pe_pp} gives the
	result. In the case (ee), combining Remark~\ref{rem:reduced_maximal_minors_zero_eq_mu_hat}, Proposition~\ref{prp:gcd_maximal_minors_phat_eq_gcd_maximal_minors_prime} (iii) and
	Proposition~\ref{prp:gcd_minors_prime_union_minors_prime_prime_succ} (ii),
	we get
	\begin{align*}
		\gcd(M^{\red}(\hat P)) & = \gcd\left(\{|\hat \mu|\} \cup \bigcup_{k=1}^{r+1}
		M^{\red}_k(\hat P)\right)                                                    \\
		                       & = \gcd\left(\{|\hat \mu|\} \cup \bigcup_{k=1}^r
		M'_k(\hat P)\right)                                                          \\
		                       & = \gcd\left(\{|\hat \mu|\} \cup M''_1(\hat
		P)\right).
	\end{align*}
	By definition, we have \( \{|\hat \mu|\} \cup M''_1(\hat P) = M'(P) \), hence
	we arrive at the claim.
\end{proof}

\begin{proof}[Proof of Theorem \ref{thm:picard_index_formula}]

	By Theorem \ref{thm:cstar_surface_from_data_completeness}, we can
	assume \( X = X(P) \subseteq Z \) as in Construction
	\ref{cns:cstar_surface_from_data}. Note that \( |\Cl(X,x)| \neq 1 \)
	only holds for fixed points of the \( \KK^* \)-action, which we
	denote by \( X^{\rm fix} \). By Remark
	\ref{rem:fixed_points_cstar_surface_of_data} and Proposition
	\ref{prp:picard_index_eq_picard_index_toric_ambient} (iii), we get
	\[
		\prod_{x \in X} |\Cl(X,x)| = \prod_{x \in X^{\fix}} |\Cl(X,x)| =
		\prod_{\sigma \in \Sigma_{\max}} |\Cl(Z,z_{\sigma})|.
	\]
	Combining Proposition
	\ref{prp:picard_index_eq_picard_index_toric_ambient} (iv) and
	Proposition \ref{prp:picard_index_of_toric_variety}, we get
	\[
		\iota_{\Pic}(X) = \iota_{\Pic}(Z) = \frac{1}{|\hat K|} \prod_{\sigma \in \Sigma_{\max}}
		|\Cl(Z,z_{\sigma})|,
	\]
	By Proposition \ref{prp:picard_index_eq_picard_index_toric_ambient}
	(i), it remains to show that \( |\hat K| = |\Cl(Z)^{\tors}| \). Since
	\( \Cl(Z) \) is isomorphic to the cokernel of \( P^* \), we have
	\[
		|\Cl(Z)^{\tors}| = \gcd(M(P^*)) = \gcd(M(P)),
	\]
	where \( M(A) \) denotes the set of maximal minors of a matrix \( A
	\). On the other hand,
	Proposition~\ref{prp:kernel_pi_eq_picard_group} along with
    Remark~\ref{rem:picard_group_free_of_cstar_surface} says that \( \hat K
	\) is the cokernel of \( \hat P^* \), hence
	\[
		|\hat K| = \gcd(M(\hat P^*)) = \gcd(M(\hat P)).
	\]
	Proposition \ref{prp:gcd_maximal_minors_P_eq_gcd_maximal_minors_Phat}
	now implies the claim.
\end{proof}

\begin{example}
	\label{exm:running_example_picard_index}
	Recall that the \( \KK^* \)-surface \( X = X(P) \) from
	Example~\ref{exm:running_example_toric_embedding} has three fixed
	points \( x^+,x^- \) and \( x_{01} \). By
	Proposition~\ref{prp:picard_index_eq_picard_index_toric_ambient}
	(iii), we can compute their local class groups via toric geometry:
	\[
		\Cl(X,x^+)\ =\ \ZZ/20\ZZ, \qquad \Cl(X,x^-)\ =\ \ZZ/12\ZZ, \qquad
		\Cl(X,x_{01})\ =\ \{0\}.
	\]
	Recall that \( \Cl(X) = \ZZ\times\ZZ/4\ZZ \). By
	Theorem~\ref{thm:picard_index_formula}, we again arrive at \( \iota_{\Pic}(X) =
	60 \).
\end{example}

Specializing Theorem~\ref{thm:picard_index_formula} to toric
surfaces, we get the following:

\begin{corollary}
	\label{cor:picard_index_formula_toric}

	Let \( Z = Z_{\Sigma} \) be a projective toric surface. Then
	\[
		\iota_{\Pic}(Z) = \frac{1}{|\Cl(Z)^{\tors}|}\prod_{\sigma \in
			\Sigma_{\max}} |\Cl(Z,z_{\sigma})|.
	\]
\end{corollary}

\begin{remark}
	\label{rem:picard_index_of_weighted_projective_plane}

	Let \( Z = \PP(w_0, w_1, w_2) \) be a weighted projective plane. We
	can assume the weights \( w_i \) to be pairwise coprime. The divisor
	class group of a weighted projective space is torsion-free and the
	orders of the local class groups at the toric fixed points are equal
	to the weights \( w_i \). By Corollary
	\ref{cor:picard_index_formula_toric}, we obtain
	\[
		\iota_{\Pic}(Z) = w_0 w_1 w_2.
	\]
	In Proposition \ref{prp:picard_index_of_fwpp}, we will generalize
	this formula to fake weighted projective planes. For weighted
	projective planes, there is also a direct way to compute the Picard
	index: The subgroup of divisor classes that are principle on the \( i
	\)-th standard affine chart \( U_i = Z \backslash V(x_i) \) is
	generated by \( w_i = [V(x_i)] \in \Cl(Z) \iso \ZZ \). For the Picard
	group, this means
	\[
		\Pic(Z)\ =\ \ZZ w_0 \cap \ZZ w_1 \cap \ZZ w_2\ =\ \ZZ
		\lcm(w_0, w_1, w_2)\ \subseteq\ \ZZ\ \iso\ \Cl(Z).
	\]
	Since the weights are pairwise coprime, we have \( \lcm(w_0, w_1,
	w_2) = w_0 w_1 w_2 \).

\end{remark}

The following example shows that Corollary
\ref{cor:picard_index_formula_toric} does not hold for higher
dimensional toric varieties.

\begin{example}
	\label{exm:toric_threefold_counter_example}
	Consider the three-dimensional weighted projective space \( Z =
	\PP(2,2,3,5) \). Note that the weights are well-formed, i.e. any three
	weights have no common factor. The Picard group is given by
	\[
		\Pic(Z)\ =\ 2 \ZZ\ \cap\ 2 \ZZ\ \cap\ 3 \ZZ\ \cap\ 5 \ZZ\ =\ 30 \ZZ \ \subseteq
		\ \ZZ \ \iso\ \Cl(Z).
	\]
	Hence we have \( \iota_{\Pic}(Z) = 30 \). On the other hand, the
	product of the orders of the local class groups is \( 60 \).
\end{example}

We now consider the \( D_8 \)-singular log del Pezzo surface of
Picard number one, which does not admit a \( \KK^* \)-action. Using
the description of its Cox Ring \cite{HaKeLa}*{Theorem 4.1}, we
construct the surface via its canonical ambient toric variety, see
also \cite{ArDeHaLa}*{Sections 3.2 and 3.3}.

\begin{example}
	\label{exm:d8_log_del_pezzo}

	Consider the integral matrix
	\[
		P\ :=\
		\begin{bmatrix}
			v_1 & v_2 & v_3 & v_4
		\end{bmatrix}\ :=\
		\begin{bmatrix}
			1 & 0 & 1 & -3 \\
			0 & 1 & 1 & -2 \\
			0 & 0 & 2 & -2
		\end{bmatrix}.
	\]
	Let \( Z = Z_{\Sigma} \) be the toric variety whose fan \( \Sigma \)
	has the following maximal cones:
	\[
		\sigma_{12}\ :=\ \cone(v_1, v_2), \qquad \sigma_{23}\ :=\ \cone(v_2,
		v_3), \qquad \sigma_{24}\ :=\ \cone(v_2,v_4),
	\]
	\[
		\sigma_{134}\ :=\ \cone(v_1,v_3,v_4).
	\]
	Let \( p \colon \hat Z \to Z \) be Cox's quotient presentation of \(
	Z \), where \( \hat Z \subseteq \bar Z := \KK^4 \). Consider the
	polynomial
	\[
		f := T_1^2 - T_2 T_3 T_4^2 + T_3^4 + T_4^4.
	\]
	We obtain a commutative diagram
	\[
		\begin{tikzcd}
			V(f) \ar[r, phantom, "=:"]                                       &
			\bar X \ar[r, hookrightarrow]                                    & \bar Z \ar[r, phantom, "="]       & \KK^4
			\\
			\bar X \cap \hat Z \ar[r, phantom, "=:"]                         & \hat X \ar[r,
			hookrightarrow] \ar[u, phantom, sloped, "\subseteq"] \ar[d, "p"] & \hat Z \ar[u,
			phantom, sloped, "\subseteq"] \ar[d, "p"]                                                                    \\
			p(\hat X) \ar[r, phantom, "=:"]                                  & X \ar[r, hookrightarrow, "\iota"]
			                                                                 & Z.
		\end{tikzcd}
	\]
	Here, \( \iota \colon X \to Z \) is the canonical toric embedding in
	the sense of \cite{ArDeHaLa}*{Sec. 3.2.5}. This implies that \( X \)
	has non-empty intersection with the toric orbits \( \TT^3 \cdot
	z_{\sigma} \) for \( \sigma \in \Sigma_{\max} \). We obtain a
	decomposition into pairwise disjoint pieces
	\[
		X = \bigcup_{\sigma \in \Sigma_{\max}} X(\sigma), \qquad \text{where}\
		X(\sigma) := X \cap \TT^3 \cdot z_{\sigma}.
	\]
	By \cite{ArDeHaLa}*{Proposition 3.3.1.5}, we have \( \Cl(X,x) \iso
	\Cl(Z,z_{\sigma}) \) for \( x \in X(\sigma) \). Note that \(
	\sigma_{12}, \sigma_{23} \) and \( \sigma_{24} \) are regular and
	\(|\Cl(Z, z_{\sigma_{134}})| = |\det(v_1, v_3, v_4)| = 2 \). This
	shows that
	\[
		\prod_{x \in X} |\Cl(X,x)|\ =\ \prod_{\sigma \in \Sigma_{\max}} |\Cl(z,
		z_{\sigma})|\ =\ 2.
	\]
	We turn to the Picard index. In the notation of Construction
	\ref{cns:picard_group_of_toric_variety}, we have
	\[
		N_{\sigma_{ij}}\ =\ \ZZ v_i + \ZZ v_j, \qquad \qquad N_{\sigma_{134}}\
		=\ N\ =\ \ZZ^3.
	\]
	Under the lattice bases \( \{v_i, v_j\} \) of \( N_{\sigma_{ij}} \),
	we can view \( P_{\sigma_{ij}} \) as the identity matrix and \(
	P_{\sigma_{134}} = \begin{bmatrix} v_1 & v_3 & v_4 \end{bmatrix} \). Computing
	matrix representations of the maps involved in Construction
	\ref{cns:picard_group_of_toric_variety}, we obtain
		{\small \[
				\begin{array}{rclcrcl}
					\alpha                                                   & =      &
					\begin{bmatrix}
						1 & 0 & 0 & 1 & 0 & -3 & 1 & 0 & 0 \\
						0 & 1 & 1 & 1 & 1 & -2 & 0 & 1 & 0 \\
						0 & 0 & 0 & 2 & 0 & -2 & 0 & 0 & 1
					\end{bmatrix}, & \qquad &
					\beta                                                    & =      &
					\begin{bmatrix}
						1 & 0 & 0 & 0 & 0 & 0 & 1 & 0 & 0 \\
						0 & 1 & 1 & 0 & 1 & 0 & 0 & 0 & 0 \\
						0 & 0 & 0 & 1 & 0 & 0 & 0 & 1 & 0 \\
						0 & 0 & 0 & 0 & 0 & 1 & 0 & 0 & 1
					\end{bmatrix},              \\\\

					\gamma                                                   & =      &
					\begin{bmatrix}
						-1 & 0  & 3 & -1 & 0  & 0  \\
						-1 & -1 & 2 & 0  & -1 & -1 \\
						0  & 0  & 0 & 0  & 0  & 1  \\
						1  & 0  & 0 & 0  & 0  & 0  \\
						0  & 1  & 0 & 0  & 0  & 0  \\
						0  & 0  & 1 & 0  & 0  & 0  \\
						0  & 0  & 0 & 1  & 0  & 0  \\
						0  & 0  & 0 & 0  & 1  & 0  \\
						-2 & 0  & 2 & 0  & 0  & 0
					\end{bmatrix},         & \qquad &
					\delta                                                   & =      &
					\begin{bmatrix}
						0  & 0  & -1 & 0  & 0  \\
						-1 & 0  & 0  & -1 & 0  \\
						0  & 0  & 0  & 1  & 0  \\
						0  & 0  & 0  & 0  & 1  \\
						1  & 0  & 0  & 0  & 0  \\
						0  & 1  & 0  & 0  & 0  \\
						0  & 0  & 1  & 0  & 0  \\
						0  & 0  & 0  & 0  & -1 \\
						0  & -1 & 0  & 0  & 0
					\end{bmatrix},

				\end{array}
			\]
			\[
				\hat P =
				\begin{bmatrix}
					0 & 0 & 0 & 0 & 1  \\
					1 & 0 & 0 & 0 & 0  \\
					0 & 1 & 0 & 0 & 0  \\
					0 & 3 & 1 & 0 & -1 \\
					0 & 2 & 0 & 0 & -1 \\
					0 & 0 & 0 & 1 & 0
				\end{bmatrix}.
			\]}
	We see that \( \hat P^* \), is surjective, hence its cokernel \( \hat
	K \) is trivial. Using
	Proposition~\ref{prp:picard_index_of_toric_variety}, we conclude \(
	\iota_{\Pic}(X) = \iota_{\Pic}(Z) = 2 \). On the other hand, we have
	\[
		\Cl(X)\ \iso\ \Cl(Z)\ \iso\ \ZZ^3 / \im(P^*)\ \iso\ \ZZ \oplus \ZZ / 2 \ZZ.
	\]
	This shows that the formula from
	Theorem~\ref{thm:picard_index_formula} does not hold for \( X \).

\end{example}

\goodbreak

\section{Log del Pezzo \texorpdfstring{\( \KK^* \)}{K*}-surfaces of Picard number one}
\label{sec:log_del_pezzo_classification}

\newcommand{\zz}{\mathbf{z}}

We contribute to the classification of log del Pezzo \( \KK^*
\)-surfaces of Picard number one, which in the toric case are
\emph{fake weighted projective planes}. As a special case of fake
weighted projective spaces, these have been studied by several
authors \cites{Ka,HaHaHaSp}. Our classification relies on Theorem
\ref{thm:picard_index_formula} to produce bounds for the number of
log del Pezzo \( \KK^* \)-surfaces with fixed Picard index. For a
related classification of fake weighted projective spaces by their
Gorenstein index, we refer to \cite{Bae}.

We begin with the toric case. Consider an integral \( 2 \times 3 \)
matrix
\[
	P = \begin{bmatrix} v_0 & v_1 & v_2 \end{bmatrix},
\]
where the columns are primitive vectors generating \( \QQ^3 \) as a
convex cone. Then we obtain a fake weighted projective plane \( Z =
Z(P) \) as the toric surface associated to the unique complete
lattice fan \( \Sigma \) with ray generators \( v_i \). We obtain a
splitting of the divisor class group:
\[
	\Cl(Z)\ \iso\ \ZZ\times\ZZ/n\ZZ.
\]
Under this splitting, we define \( (w_i,\eta_i) := \omega_i := [D_i]
\), where \( D_i \) is the torusinvariant Weil divisor associated to
\( v_i \), for \( i=0,1,2 \). Note that any two of the~\( \omega_i \)
generate~\( \Cl(Z) \). In particular, the \( w_i \) are pairwise
coprime integers.

\begin{proposition}
	\label{prp:picard_index_of_fwpp}

	Let \( Z = Z(P) \) be a fake weighted projective plane with divisor
	class group \( \Cl(Z) \iso \ZZ \times \ZZ/n\ZZ \). Write \( \zz(i) \)
	for the toric fixed points associated to \( \sigma_i := \cone(v_j\ ;\
	j\neq i) \). Then the orders of the local class groups are given by
	\[
		\vert \Cl(Z,\zz(i)) \vert = nw_i = \vert \det(\sigma_i) \vert.
	\]
	Moreover, the Picard index is \(\iota_{\Pic}(Z) = n^2 w_0 w_1 w_2.\) \end{proposition}

\begin{proof}
	We have \( \Cl(Z,\zz(i))\ \iso\ K / \ZZ\omega_i \),
	hence the order equals \( nw_i \). On the other hand, \(
	\Cl(Z,\zz(i)) \iso \Cl(U_i) \), where \( U_i \) is the affine toric
	chart associated to \( \sigma_i \), hence the order equals \( \vert
	\det(\sigma_i) \vert \). The formula for the Picard Index follows
	from Corollary~\ref{cor:picard_index_formula_toric}.
\end{proof}

\begin{remark}
	\label{rem:fwpp_generator_matrix}

	Let \( Z(P) \) be a fake weighted projective plane. Then we have \(
	Z(P) \iso Z(P') \) if and only if \( P' = A \cdot P \cdot S \) holds
	with a unimodular matrix \( A \) and a permutation matrix \( S \).

\end{remark}

\begin{proposition}
	\label{prp:generator_matrix_of_fwpp}

	Let \( Z(P) \) be a fake weighted projective plane. Then there exists
	an integer \( 0 \leq x < nw_2 \) with \( \gcd(x, nw_2) = 1 \) such
	that \( Z(P) \iso Z(P') \), where
	\[
		P' =
		\begin{bmatrix}
			1 & x    & - \frac{w_0 + x w_1}{w_2} \\
			0 & nw_2 & -nw_1
		\end{bmatrix}.
	\]
\end{proposition}

\begin{proof}
	We write \( P = \begin{bmatrix} v_0 & v_1 & v_2 \end{bmatrix} \),
	where \( v_i = (x_i, y_i) \in \ZZ^2 \) are primitive
	vectors. Applying a unimodular matrix from the left achieves \( v_0 = (1,0)
	\) and \( 0 \leq x_1 < y_1 \). By Proposition
	\ref{prp:picard_index_of_fwpp}, we have
	\[
		nw_2 = \det(v_0, v_1) = y_1, \qquad \qquad -nw_1 = \det(v_0, v_2) = y_2.
	\]
	\[
		nw_0 = \det(v_1, v_2) = -nw_1x_1 - nw_2x_2.
	\]
	This shows the claim.
\end{proof}

\begin{algorithm}
	\label{alg:classify_fwpp_by_picard_index}

	\emph{Input:} A positive integer \( \iota \), the prospective Picard Index.
	\emph{Algorithm:}
	\begin{itemize}
		\item Set \( L := \emptyset \).
		\item For each quadruple \( (n,w_0,w_1,w_2) \) of positive integers with \(
		      n^2w_0w_1w_2 = \iota \) such that \( (w_0,w_1,w_2) \) are pairwise
		      coprime, do:
		      \begin{itemize}
			      \item[\(\bullet\)] For each \( 0 \leq x < nw_2 \) with \( \gcd(x,nw_2) = 1 \), do:
				      \begin{itemize}
					      \item[\(\bullet\)] Set
						      \[
							      P :=
							      \begin{bmatrix}
								      1 & x    & - \frac{w_0 + x w_1}{w_2} \\
								      0 & nw_2 & -nw_1
							      \end{bmatrix}.
						      \]
					      \item[\(\bullet\)] If for all permutation matrices \( S \),
						      the Hermite normal form of \( P \cdot S \) differs
						      from the Hermite normal form of all \( P' \in L \),
						      add \( P \) to \( L \).
				      \end{itemize}
			      \item[\(\bullet\)] end do.
		      \end{itemize}
		\item end do.
	\end{itemize}
	\emph{Output:} The set \( L \). Then every fake weighted projective plane \(
	Z\) with \( \iota_{\Pic}(Z) = \iota \) is isomorphic to precisely one \(
	Z(P) \) with \( P \in L \).
\end{algorithm}

\begin{proof}
	Propositions \ref{prp:picard_index_of_fwpp} and
	\ref{prp:generator_matrix_of_fwpp} show that the algorithm produces all
	generator matrices of fake weighted projective planes sharing \( \iota \) as
	their Picard index. By Remark \ref{rem:fwpp_generator_matrix}, the matrices
	in \( P \in L \) belong to pairwise non-isomorphic surfaces.
\end{proof}

Using this Algorithm, we arrive at the following classification.

\begin{theorem}
	\label{thm:classification_fwpp}
	There are \( 15\,086\,426 \) isomorphy classes of fake weighted projective planes
	with Picard index at most \( 1\,000\,000 \). Of those, \( 68\,053 \) have
	Picard index at most \( 10\,000 \). The number of isomorphy classes for
	given Picard index develops as follows:

	\begin{center}
		\includegraphics{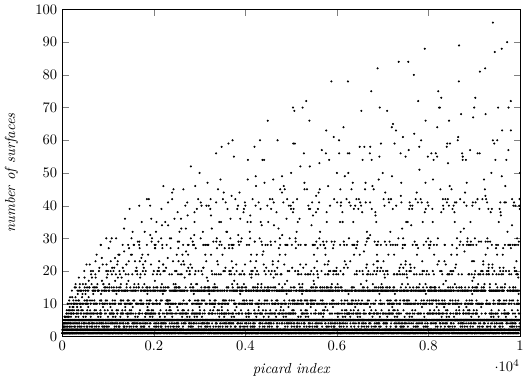}
	\end{center}

\end{theorem}

We turn to non-toric rational projective \( \KK^* \)-surfaces of
Picard number one. These arise from defining matrices as in
Construction \ref{cns:cstar_surface_from_data}, see also
\cite{HaHaHaSp}*{Section 4}.

Given a rational projective \( \KK^* \)-surface \( X = X(P) \subseteq
Z \) as in Construction~\ref{cns:cstar_surface_from_data}, the
columns \( v_{ij},v^{\pm} \) define torusinvariant Weil divisors \(
D_{ij},D^{\pm} \) on \( Z \). If \( X \) is of Picard number one, the
divisor class group splits as
\[
	\Cl(X)\ \iso\ \Cl(Z)\ \iso\ \ZZ \times \Cl(Z)^{\tors}.
\]
We define \( (w_{ij},\eta_{ij}) := \omega_{ij} := [D_{ij}] \) and \(
(w^{\pm},\eta^{\pm}) := \omega^{\pm} := [D^{\pm}] \).

\begin{proposition}
	\label{prp:picard_index_of_cstar_surface_rho_one_ee}

	Let \( X = X(P) \), where \( P = \begin{bmatrix} v_{01} & v_{02} &
                v_1    & \dots  & v_r\end{bmatrix} \) is a defining matrix of type
	(ee) with \( n_0 = 2 \) and \( n_1 = \dots = n_r = 1 \). Then \( X \)
	is of Picard number one. We have two elliptic
	fixed points \( x^{\pm} \) and one hyperbolic fixed
	point~\( x_{01} \). With \( \lambda := |\Cl(X)^{\tors}| \), the orders of
	the local class groups are
	\[
		|\Cl(X,x^+)|\ =\ \lambda w_{02} = |\det(\sigma^+)|, \qquad
		|\Cl(X,x^-)|\ =\ \lambda w_{01} = |\det(\sigma^-)|.
	\]
	\[
		|\Cl(X,x_{01})| = \det
		\begin{bmatrix}
			-l_{01} & -l_{02} \\
			d_{01}  & d_{02}
		\end{bmatrix} = d_{01} l_{02} - d_{02} l_{01}.
	\]
	With \( M := |\Cl(X,x_{01})| \), the Picard index is
	\(\iota_{\Pic}(X) = \lambda w_{01} w_{02} M.\)

\end{proposition}

\begin{proposition}
	\label{prp:picard_index_of_cstar_surface_rho_one_ep}

	Let \( X = X(P) \), where \( P = \begin{bmatrix} v_0 & \dots & v_r & v^- \end{bmatrix} \) is a
	defining matrix of type (ep) with \( n_0 = \dots = n_r = 1 \).
	Then \( X = X(P) \) is of
	Picard number one. We have one elliptic fixed point \( x^+ \) and
	parabolic fixed points \( x_{in_i} \) for \( i = 0, \dots, r \). With
	\( \lambda := |\Cl(X)^{\tors}| \), the orders of the local class
	groups are
	\[
		|\Cl(X,x^+)| = \lambda w^- = \det(\sigma^+), \qquad
		|\Cl(X,x_{in_i})| = \det
		\begin{bmatrix}
			l_i & 0 \\
			d_i & 1
		\end{bmatrix} = l_i
	\]
	Moreover, the Picard index is \(\iota_{\Pic}(X) = w^- l_0 \cdots
	l_r.\)

\end{proposition}

\begin{proof}[Proof of Propositions
		\ref{prp:picard_index_of_cstar_surface_rho_one_ee} and
		\ref{prp:picard_index_of_cstar_surface_rho_one_ep}]

	In both cases, we have \( n+m = r+2 \), hence~\( X(P) \) is of Picard
	number one. For the descriptions of the fixed points, see Remark
	\ref{rem:fixed_points_cstar_surface_of_data}. For the orders of the
	local class groups, we apply \cite{ArDeHaLa}*{Lemma 2.1.4.1} to the
	toric ambient variety of \( X \) and use
	Proposition~\ref{prp:picard_index_eq_picard_index_toric_ambient}
	(iii). Applying Theorem \ref{thm:picard_index_formula} gives the
	expressions for the Picard index.
\end{proof}

\begin{proposition}[See {\cite{HaHaHaSp}*{Propositions 5.9 and 5.10}}]
	\label{prp:platonic_triples_of_log_terminality}
	Let \( X \) be a non-toric rational projective \( \KK^* \)-surface of Picard number
	one. Assume that \( X \) is log terminal. Then we have \( X \iso X(P) \)
	for some defining matrix \( P \) as in one of the Propositions
	\ref{prp:picard_index_of_cstar_surface_rho_one_ee} and
	\ref{prp:picard_index_of_cstar_surface_rho_one_ep}. The possible
	tuples \( (l_{01}, l_{02}, l_1, \dots, l_r) \) for type (ee) and \( (l_0,
	\dots, l_r) \) for type (ep) are the following:
	\[
		\begin{array}{cllcll}
			\mathrm{(eAeA):} & (1,1,x_1,x_2), &                & \mathrm{(eAeD):} &
			(1,y,2,2),       & (1,2,y,2),
			\\[5pt]
			\mathrm{(eAeE):} & (1,z,3,2),     & (1,3,z,2),     &
			\mathrm{(eDeD):} & (2,2,y,2),     & (y_1,y_2,2,2)
			\\
			                 & (1,2,z,3),     &                &
			                 & (1,1,y,2,2),   &                                              \\[5pt]
			\mathrm{(eDeE):} & (2,3,z,2),     & (2,z,3,2)      &
			\mathrm{(eEeE):} & (2,2,z,3)      & (z_1,z_2,3,2),                               \\
			                 &                &                &
			                 & (3,3,z,2)      & (1,1,z,3,2),                                 \\[5pt]
			\mathrm{(eDp):}  & (y,2,2),       &                & \mathrm{(eEp):}  & (z,3,2),
			                 &                                                               \\
		\end{array}
	\]
	where \( 2 \leq x_1, x_2, y, y_1, y_2 \) and \( 3 \leq z, z_1 z_2
	\leq 5 \) and the notation ``eA, eD, eE'' refers to log terminal \(
	x^{\pm} \in X \) of type \( A_n, D_n \) or \( E_6, E_7, E_8 \).
	Moreover, any non-toric, rational, log terminal, projective \( \KK^*
	\)-surface of Picard number one is del Pezzo.
\end{proposition}

Using the formula for the Picard index in Proposition
\ref{prp:picard_index_of_cstar_surface_rho_one_ee}, we can go through
all possible quadruples of positive integers \( (\lambda, w_{01},
w_{02}, M) \) whose product equals a given Picard index. Going
through all the possible tuples \( (l_{01}, l_{02}, l_1, \dots, l_r)
\) in Proposition \ref{prp:platonic_triples_of_log_terminality}, we
can use the equations for the orders of the local class groups as
well as the fact that \( (w_{01},w_{02},w_{11},\dots,w_{r1}) \in
\ker(P) \) to derive bounds for all remaining entries of \( P \).
This yields an efficient classification of all non-toric rational
projective log terminal \( \KK^* \)-surfaces of Picard number one and
type (ee) by their Picard index. The type (ep) is handled
analogously. We arrive at the following classification.

\begin{theorem}
	\label{thm:classification_nontoric_cstar}
	There are \( 1\,347\,433 \) families of non-toric, log del Pezzo
	\( \KK^* \)-surfaces of Picard number one and Picard index at
	most \( 10\,000 \). The numbers of families for given Picard index develop as
	follows:

    \begin{center}
		\includegraphics{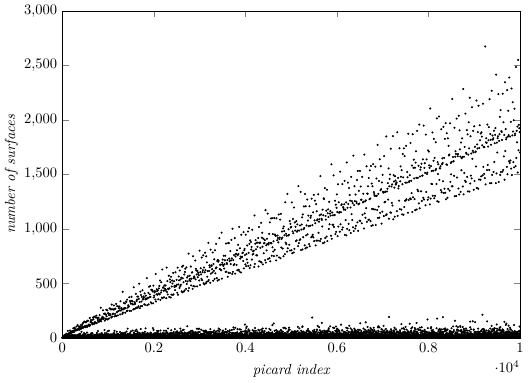}
	\end{center}
	
\end{theorem}

\begin{remark}
	\label{rem:shape_of_graph_nontoric}

Comparing the numbers of surfaces per given Picard index in the
different subcases of our classification, we observe that the
non-toric \( \KK^* \)-surfaces quickly outnumber the toric ones. The
fastest growing subcase seems to be ``eDeD'', which is responsible
for the cone shape in the upper part of the plot in
Theorem~\ref{thm:classification_nontoric_cstar}. The following table
gives an impression of the proportion of the different cases:

	\begin{center}

		{\small
			\begin{tabular}{c|ccccccccc}
				\( \iota_{\Pic} \) & toric   & eAeA     & eAeD    & eAeE   & eDeD        & eDeE & eEeE &
				eDp                & eEp                                                                        \\\hline
				\( \leq 10 \)      & 14      & 5        & 4       & 10     & 1           & 0    & 0    & 1  & 0 \\
				\( \leq 100 \)     & 243     & 260      & 129     & 39     & 117         & 4    & 15   & 28 & 5
				\\
				\( \leq 1\,000 \)  & 4\,205  & 7\,425   & 2\,209  & 206    & 11\,622     & 32   & 103
				                   & 521     & 51                                                               \\
				\( \leq 10\,000 \) & 68\,053 & 157\,482 & 31\,561 & 1\,011 & 1\,148\,587
				                   & 197     & 569      & 7\,520  & 506
			\end{tabular}
		}

	\end{center}
    Up to Picard index 2\,500, the complete list of
    defining matrices is available at~\cite{ldp}.

\end{remark}


\begin{bibdiv}

	\begin{biblist}

		\bib{AlNi}{article}{
			author={Alekseev, V. A.},
			author={Nikulin, V. V.},
			title={Classification of del Pezzo surfaces with log-terminal
					singularities of index $\le 2$ and involutions on $K3$ surfaces},
			language={Russian},
			journal={Dokl. Akad. Nauk SSSR},
			volume={306},
			date={1989},
			number={3},
			pages={525--528},
			issn={0002-3264},
			translation={
					journal={Soviet Math. Dokl.},
					volume={39},
					date={1989},
					number={3},
					pages={507--511},
					issn={0197-6788},
				},
		}

		\bib{ArDeHaLa}{book}{
			author={Arzhantsev, Ivan},
			author={Derenthal, Ulrich},
			author={Hausen, J\"urgen},
			author={Laface, Antonio},
			title={Cox rings},
			series={Cambridge Studies in Advanced Mathematics},
			volume={144},
			publisher={Cambridge University Press, Cambridge},
			date={2015},
			pages={viii+530},
		}

        \bib{Bae}{arXiv}{
              title={Sharp volume and multiplicity bounds for Fano simplices}, 
              author={Andreas Bäuerle},
              year={2023},
              eprint={2308.12719},
              archivePrefix={arXiv},
              primaryClass={math.CO},
        }

		\bib{CoLiSc}{book}{
			author={Cox, David A.},
			author={Little, John B.},
			author={Schenck, Henry K.},
			title={Toric varieties},
			series={Graduate Studies in Mathematics},
			volume={124},
			publisher={American Mathematical Society, Providence, RI},
			date={2011},
			pages={xxiv+841},
		}

		\bib{FuYa}{article}{
			author={Fujita, Kento},
			author={Yasutake, Kazunori},
			title={Classification of log del Pezzo surfaces of index three},
			journal={J. Math. Soc. Japan},
			volume={69},
			date={2017},
			number={1},
			pages={163--225},
			issn={0025-5645},
		}

		\bib{Ful}{book}{
			author={Fulton, William},
			title={Introduction to toric varieties},
			series={Annals of Mathematics Studies},
			volume={131},
			note={The William H. Roever Lectures in Geometry},
			publisher={Princeton University Press, Princeton, NJ},
			date={1993},
			pages={xii+157},
			isbn={0-691-00049-2},
		}

		\bib{Hae}{book}{
		author={H\"{a}ttig, Daniel},
		title={Lattice polygons and surfaces with torus action},
		series={PhD Thesis},
		publisher={Universit\"at T\"ubingen},
		date={2023},
		}

		\bib{HaHaHaSp}{arXiv}{
		author={H\"{a}ttig, Daniel},
		author={Hafner, Beatrice},
		author={Hausen, J\"{u}rgen},
		author={Justus Springer},
		title={Del Pezzo surfaces of Picard number one admitting a torus action},
		year={2022},
		eprint={2207.14790},
		archivePrefix={arXiv},
		primaryClass={math.AG}
		}

		\bib{ldp}{webpage}{
		title={ldp-database},
		subtitle={Log del Pezzo surfaces with torus action - a searchable
				database},
		date={2023},
		accessdate={October 2023},
		url={https://www.math.uni-tuebingen.de/forschung/algebra/ldp-database/},
		author={Daniel H\"{a}ttig},
		author={Hausen, J\"{u}rgen},
		author={Justus Springer}
		}

		\bib{HaHaeSp}{arXiv}{
		title={Classifying log del Pezzo surfaces with torus action},
		author={Daniel H\"{a}ttig},
		author={Hausen, J\"{u}rgen},
		author={Justus Springer},
		year={2023},
		eprint={2302.03095},
		archivePrefix={arXiv},
		primaryClass={math.AG}
		}

		\bib{HaHe}{article}{
		author={Hausen, J\"{u}rgen},
		author={Herppich, Elaine},
		title={Factorially graded rings of complexity one},
		conference={
		title={Torsors, \'{e}tale homotopy and applications to rational points},
		},
		book={
				series={London Math. Soc. Lecture Note Ser.},
				volume={405},
				publisher={Cambridge Univ. Press, Cambridge},
			},
		date={2013},
		pages={414--428},
		}

		\bib{HaHeSu}{article}{
		author={Hausen, J\"{u}rgen},
		author={Herppich, Elaine},
		author={S\"{u}\ss, Hendrik},
		title={Multigraded Factorial Rings and Fano Varieties with Torus Action},
		journal={Documenta Math.},
		volume={16},
		date={2011},
		number={3},
		pages={71--109},
		issn={1431-0635},
		}

		\bib{HaSu}{article}{
		author={Hausen, J\"{u}rgen},
		author={S\"{u}\ss , Hendrik},
		title={The Cox ring of an algebraic variety with torus action},
		journal={Adv. Math.},
		volume={225},
		date={2010},
		number={2},
		pages={977--1012},
		issn={0001-8708},
		}

		\bib{HaKeLa}{article}{
			ISSN = {00255718, 10886842},
			URL = {https://www.jstor.org/stable/mathcomp.85.297.467},
			author = {Jürgen Hausen and Simon Keicher and Antonio Laface},
			journal = {Mathematics of Computation},
			number = {297},
			pages = {467--502},
			publisher = {American Mathematical Society},
			title = {Computing Cox rings},
			volume = {85},
			year = {2016}
		}

		\bib{HaWr}{article}{
		author={Hausen, J\"{u}rgen},
		author={Wrobel, Milena},
		title={Non-complete rational $T$-varieties of complexity one},
		journal={Math. Nachr.},
		volume={290},
		date={2017},
		number={5-6},
		pages={815--826},
		issn={0025-584X},
		}

		\bib{Ka}{article}{
			author={Kasprzyk, Alexander M.},
			title={Bounds on fake weighted projective space},
			journal={Kodai Math. J.},
			volume={32},
			date={2009},
			number={2},
			pages={197--208},
			issn={0386-5991},
			review={\MR{2549542}},
		}

		\bib{Nak}{article}{
			author={Nakayama, Noboru},
			title={Classification of log del Pezzo surfaces of index two},
			journal={J. Math. Sci. Univ. Tokyo},
			volume={14},
			date={2007},
			number={3},
			pages={293--498},
			issn={1340-5705},
		}

		\bib{OrWa1}{article}{
			author={Orlik, Peter},
			author={Wagreich, Philip},
			title={Isolated singularities of algebraic surfaces with $\mathbb{C}^*$-action},
			journal={Ann. of Math. (2)},
			volume={93},
			date={1971},
			pages={205--228},
			issn={0003-486X},
		}

		\bib{OrWa2}{article}{
			author={Orlik, Peter},
			author={Wagreich, Philip},
			title={Singularities of algebraic surfaces with $\mathbb{C}^*$-action},
			journal={Math. Ann.},
			volume={193},
			date={1971},
			pages={121--135},
			issn={0025-5831},
		}

		\bib{OrWa3}{article}{
			author={Orlik, Peter},
			author={Wagreich, Philip},
			title={Algebraic surfaces with $k^*$-action},
			journal={Acta Math.},
			volume={138},
			date={1977},
			number={1-2},
			pages={43--81},
			issn={0001-5962},
		}

	\end{biblist}

\end{bibdiv}
\end{document}